\documentclass[12pt]{amsart}
\usepackage{amsmath}

\theoremstyle{plain}
\newtheorem{theorem}{Theorem}[section]
\newtheorem{lemma}[theorem]{Lemma}
\newtheorem{corollary}[theorem]{Corollary}
\newtheorem{proposition}[theorem]{Proposition}
\newtheorem{ther}{Theorem}

\theoremstyle{remark}
\newtheorem*{remark}{Remark}

\newtheorem*{ack}{Acknowledgements}

\theoremstyle{definition}
\newtheorem{defn}[theorem]{Definition}

\def \R {\mathbf{R}}
\def \Z {\mathbf{Z}}

\def \i {\mathbf{i}}

\def\cp{\mathbf{CP}}
\def\cpbar{\overline{\mathbf{CP}}^{2}}

\def\AA{\mathcal{A}}
\def\HH{\mathcal{H}}

\def\GG{\mathcal{G}}

\def\RR{\mathcal{R}}

\def\MM{\mathcal{M}}
\def\NN{\mathcal{N}}
\def\DD{\mathcal{D}}
\def\FF{\mathcal{F}}
\def\VV{\mathcal{V}}
\def\SS{\mathcal{S}}

\def\delN{\partial N}

\def\hatN{\widehat N}
\def\dtau{\mathbf{\dot{\mathrm{\tau}}}}

\DeclareMathOperator\im{im}
\DeclareMathOperator\sign{sign}

\DeclareMathOperator\ind{ind}
\DeclareMathOperator\Spin{Spin}
\def\spin{\ifmmode{\Spin}\else{$\Spin$\ }\fi}
\DeclareMathOperator{\Spinc}{Spin^c}
\def\spinc{\ifmmode{\Spinc}\else{$\Spinc$\ }\fi}
\def\sw{Seiberg-Witten}
\def\SW{\ifmmode{\text{SW}}\else{$\text{SW}$}\fi}
\def\swi{Seiberg-Witten invariant}
\def\swe{Seiberg-Witten equations}
\def\sws{Seiberg-Witten solutions}

\def\iso{\cong}
\def\col{\colon}

\def\vfi{\varphi}
\def\veps{\varepsilon}

\def\alfahat{\widehat\alpha}
\def\etahat{\widehat\eta}

\newlength{\zamik}
\newlength{\sirodst}
\def\odst#1#2{\hbox to\zamik{(#1)}\begin{minipage}[t]{\sirodst}#2\end{minipage}\\[4mm]}

\def\nad#1{\stackrel{#1}{\longrightarrow}}
\def\der#1#2{\frac{\partial #1}{\partial #2}}
\def\incl{\hookrightarrow}
\def\sto{\kern-2mm{\to}\kern-2mm}
\def\snad#1{\kern-2mm\stackrel{#1}{\longrightarrow}\kern-2mm}

\def\dsp{\displaystyle}

\begin{document}

\title{Bounds on genus and geometric intersections from cylindrical end moduli spaces}

\author[Sa\v{s}o Strle]{Sa\v{s}o Strle}
\address{Department of Mathematics and Statistics, McMaster University, Hamilton, ON L8S 4K1, Canada}
\email{\rm{strles@math.mcmaster.ca}}
\thanks{A part of this work was done while the author was a Liftoff mathematician with Clay Mathematics Institute.
The author was partially supported by MSZS of the Republic of Slovenia Research Program No.~101-509.}
\date\today
\begin{abstract}
In this paper we present a way of computing a lower bound for genus of any smooth representative of a
homology class of positive self-intersection in a smooth four-manifold $X$ with second positive Betti number
$b_2^+(X)=1$. We study the solutions of the \swe\ on the cylindrical end manifold which is the complement of the
surface representing the class. The result can be formulated as a form of generalized adjunction inequality.
The bounds obtained depend only on the rational homology type of the manifold, and include the Thom conjecture
as a special case. We generalize this approach to derive lower bounds on the number of intersection points of $n$
algebraically disjoint surfaces of positive self-intersection in manifolds with $b_2^+(X)=n$.
\end{abstract}

\maketitle

\markboth{SA\v{S}O STRLE}{BOUNDS ON GENUS AND GEOMETRIC INTERSECTIONS}


\section*{Introduction}
\sw\ theory has proved very useful in the study of the minimal genus problem.
After Kronhei\-mer-Mrowka's proof of the Thom conjecture \cite{kronheimer-mrowka:thom},
regarding the minimal genus problem in the complex projective plane, the following result (the so-called generalized
Thom conjecture) was obtained by Morgan-Szabo-Taubes \cite{mst} (for classes of non-negative self-intersection)
and Ozsvath-Szabo \cite{ozsvath-szabo:thom} (the general case): any smooth symplectic curve in a closed symplectic
four-mani\-fold minimizes the genus
in its homology class. Their proofs depend on results of Taubes \cite{taubes:symplectic-sw} about
Seiberg-Witten theory of symplectic manifolds, specifically the basic classes of such manifolds. In this paper
we present a way of deriving genus bounds that does not depend on any special structure on the manifold.
Rather than working over a closed manifold, we study \swe\ on an associated cylindrical end manifold; this approach
is related to the work of Fr{\o}yshov \cite{froyshov}. The bounds obtained in this way
depend only on the rational homology type of the manifold. Since the results are so general, the information
about a possible symplectic structure on $X$ is lost; in particular, the bounds are independent of
(the sign of) the canonical class of $X$, whereas the bounds coming from the generalized Thom conjecture
detect differences in canonical classes. An important advantage of our approach is that it can be used in
manifolds with vanishing \swi{s}, in particular to study geometric intersections of surfaces.

Consider a divisible homology class $d\xi \in H_2(X;\Z)$ of positive self-inter\-sect\-ion in a smooth
four-manifold $X$ with $b_1(X)=0$ and $b_2^+(X)=1$. The divisibility $d>1$ of the homology class is
crucial (for technical reasons) while studying \swe\ on the cylindrical end manifold $Z=X-\Sigma$,
where $\Sigma$ is a smooth embedded surface representing $d\xi$. The end of $Z$ is modeled on a non-trivial
circle bundle $Y$ over $\Sigma$, and we work with \sws\ on $Z$ that exponentially decay to solutions on $Y$.
This depends on description of the perturbed \sw\ moduli spaces on $Y$ obtained by Mrowka-Ozsvath-Yu \cite{moy}.

Even though the method requires us to consider a divisible class, the main result concerning genus bounds holds
for primitive classes as well. The bound can be stated in the form of generalized adjunction inequality.
\begin{ther}\label{adjunction-inequality}
Let $X$ be a smooth closed oriented four-manifold with $b_1(X)\allowbreak =0$ and $b_2^+(X)=1$. If $\Sigma \subset X$ is a
smooth embedded surface of positive self-intersection, then
\begin{equation}\label{adjunction}
\chi(\Sigma)+[\Sigma]^2 \le |\langle c,[\Sigma]\rangle|
\end{equation}
for any characteristic vector $c \in H^2(X)$ that satisfies $c^2 > \sigma(X)$.
\end{ther}
Based on this inequality it is straightforward to derive minimal genus formulae in (rational homology)
$\cp^2$, $S^2 \times S^2$ and $\cp^2 \# \cpbar$. In manifolds with rational homology of rational surfaces
$\cp^2\# n\cpbar$ with $2 \le n \le 9$ the results are easiest to state for reduced classes, defined by Li-Li
\cite{li:mingenus}. In particular, we prove that for any $g>0$ there are only finitely many reduced classes of
minimal genus $g$ (see Proposition \ref{neg-signature}; also see \cite{li:mingenreps}). We note that genuine rational
surfaces mentioned above are `genus-minimal' in the sense that minimal genus representatives in these manifolds
have the smallest possible genus among all manifolds with the same rational homology type.

Above considerations generalize to manifolds $X$ with $b_2^+(X)=n$ in a way that allows us to study a collection
of $n$ surfaces in $X$. The counterpart of the adjunction inequality is the following result.
\begin{ther}\label{n-adj-ther}
Let $X$ be a smooth closed connected four-manifold with $b_1(X)\allowbreak =0$ and $b_2^+(X)=n>1$, and let  $\Sigma_1,
\ldots,\Sigma_n$ be disjoint embedded surfaces in $X$
with positive self-intersections. If $c \in H^2(X)$ is a characteristic vector satisfying
$$c^2>\sigma(X)\ \ \hbox{and}\ \ \langle c,[\Sigma_i] \rangle \ge 0 \ \ \hbox{for all $i$},$$
then
\begin{equation}\label{n-adjunction}
\chi(\Sigma_i)+[\Sigma_i]^2 \le \langle c,[\Sigma_i]\rangle
\end{equation}
holds for at least one $i$.
\end{ther}
We use this to derive a lower bound on the number of intersection points of surfaces of low genus. For example,
suppose that classes $(p,q,0,0)$ and $(0,0,r,s)$ in $H_2(S^2 \times S^2 \#  S^2 \times S^2)$ are represented by
spheres in the connected sum $S^2 \times S^2 \#  S^2 \times S^2$. If $p,q,r,s\ge 2$ and $p+q\ge r+s$, then the
number of intersection points of the two spheres is at least
$$pq+(r-1)(s-1).$$
In general, the lower bound on the number of intersection points obtained in this way is roughly by a factor
of 2 better than the bounds obtained via the $g$-signature Theorem (see \cite{gilmer:thesis}).
We also give an example where the bound on the number of intersection points is optimal.

This paper is divided in two parts. Part I is concerned with technical aspects of \sw\ moduli spaces
over  cylindrical end manifolds.  The main results of this part are the dimension formula for the moduli
space (Corollary \ref{formal-dim}) and compactness and regularity results of Section \ref{comp-reg}.
In Part II we use results of Part I to derive genus bounds and bounds on the number of intersection
points of surfaces. We first present a derivation of a genus bound in $\cp^2$ (which is equivalent to the Thom
conjecture) and then proceed to the
general case. This is described in Theorem \ref{main-theorem}, which can be rephrased as a generalized adjunction
inequality stated above (see Section \ref{constr}). After that we consider several examples in which one can
derive explicit formulae for genus bounds and address the question of representability. In the last
section of the paper we turn to geometric intersections of surfaces.

\begin{ack}
I would like to thank my advisor Danny Ruberman for the numerous conversations and for all his support and guidance.
I am grateful to Tom Mrowka and Jerry Levine for suggestions that helped in the preparation of this work.
I would also like to thank Professor Jo\v{z}e Vrabec for all that he taught me and for his continued support.
\end{ack}


\vskip5mm
\centerline{\large \sc Part I: Cylindrical end moduli spaces}
\section{The setup}\label{problem}
Throughout $X$ will denote a smooth closed connected oriented four-mani\-fold. If a smooth oriented surface $\Sigma$ is
embedded in $X$ so that the image of its fundamental homology class $[\Sigma]$ in $H_2(X)$ is not a torsion class, we
say that $\Sigma$ {\em represents} this homology class.
Denote by $N \subset X$ a compact tubular neighborhood of $\Sigma$.
It will be convenient to distinguish between $\delN$, oriented as the boundary of $N$, and $Y=\overline\delN$,
oriented as the `boundary' of $Z$. More precisely, let $Z_0$ be the closure of the complement of $N$ in $X$; then
$Y=\partial Z_0$ and we think of $Z=X-\Sigma$ as $Z_0$ with a half-infinite cylinder attached, $Z=Z_0 \cup_Y [0,\infty)
\times Y$. We refer to $[0,\infty) \times Y$ as the {\em cylindrical end} of $Z$ and say that the end of $Z$ is
{\em modeled on} $Y$. The following proposition summarizes the relevant cohomological information about these spaces.
Unless specified otherwise all the (co)homology groups have integer coefficients.

\begin{proposition}\label{cohomology-groups}
Let $X$ be a closed oriented four-manifold and let $\Sigma \subset X$ be an embedded surface representing the class
$d\xi$, where $\xi \in H_2(X)$ is a primitive class of non-zero self-intersection and $d\ge1$ is an integer. Denote by
$N$ a compact tubular neighborhood of $\Sigma$ and by $n=(d\xi)^2$ the degree of the circle bundle $\delN \to \Sigma$.
Then
$$H^1(\delN)\iso H^1(\Sigma), \ \ H^2(\delN) \iso H^2(\Sigma)/n[\Sigma]^* \oplus H^1(\Sigma),$$
$$H^1(Z)\iso H^1(X), \ \ H^2(Z)\iso H^2(X)/d\alpha \oplus F,$$
where $[\Sigma]^*$ denotes the fundamental cohomology class of $\Sigma$, $\alpha$ is the Poincar\'e dual of $\xi$
and  $F$ is a subgroup of $H^1(\Sigma)$. The restriction homomorphism $H^2(X) \to H^2(N)$ sends $\alpha$ to
$d\xi^2[\Sigma]^*$ and its image is a subgroup of index $d$. Moreover, the restriction homomorphism $H^2(Z) \to
H^2(\delN)$  is injective on $F$.
\end{proposition}
\begin{proof}
The cohomology groups of $\delN$ follow easily from the Gysin exact sequence of the circle  bundle
$S^1 \incl \delN \nad{p}  \Sigma$ with the Chern class $c_1=n[\Sigma]^*$.

To determine the cohomology of $Z$, use the Poincar\'e duality and excision isomorphisms:
$H^2(Z)\iso H_2(Z,\delN) \iso H_2(X,N)$. The last group can be computed using the exact sequence of the pair
$(X,N)$:
$$H_3(X) \nad{\iso} H_3(X,N) \to H_2(N) \to H_2(X) \to H_2(X,N) \to H_1(N).$$
Since $H_1(N)$ is a free abelian group, we have $H_2(X,N) \iso H_2(X)/d\xi \oplus F'$ for some subgroup $F'$ of
$H_1(N)$. Similarly we obtain $H^1(Z)\iso H_3(X,N) \iso H^1(X)$.

From $\langle\alpha,[\Sigma]\rangle=\xi \cdot d\xi$ it follows that $\alpha \in H^2(X)$ restricts to
$d\xi^2[\Sigma]^* \in H^2(N)$. The remaining claims follow easily.
\end{proof}


\section{Seiberg-Witten solutions over a circle bundle}\label{circle-bundle}
We describe the structure of the moduli spaces of solutions  of certain  {\em perturbed} \swe\ on a circle
bundle $p \col Y \to \Sigma$ of degree $n \ne 0$ over an oriented smooth surface $\Sigma$, studied by
Mrowka-Ozsvath-Yu \cite{moy}. The purpose of the perturbation is to make the equations behave as if the bundle $Y$
were a product. This is achieved by choosing a `product' connection as the background connection in $T^*Y$ in place
of the Levi-Civita connection.

To define the background connection, choose a constant curvature metric $g_\Sigma$ of volume 1 on $\Sigma$; denote by
$vol_\Sigma$ the corresponding volume form and let $\omega= p^*(vol_\Sigma)$ be its pull-back to $Y$. The circle
bundle $Y$ admits a connection 1-form $i\vfi \col TY \to i\R$ of constant curvature; observe that $d\vfi=-2\pi n\omega$,
since $Y$ has degree $n$.  This connection determines a splitting $T^*Y=\R\vfi\oplus H$, where $H\iso p^*T^*\Sigma$ is
the horizontal distribution. A metric on $Y$, compatible with this splitting, is given by $g_Y=\vfi^2 + p^*g_\Sigma$;
the corresponding volume form is $vol_Y=\vfi\wedge\omega$. Note that the radius of a fiber circle with respect to this
metric is 1. The {\em product connection} on $Y$ is defined by $\nabla^Y = d\oplus  p^*\nabla^\Sigma$, where
$\nabla^\Sigma$ is the Levi-Civita connection of $(\Sigma,g_\Sigma)$. Connection $\nabla^Y$ is compatible with the
splitting and the metric. However, it is not torsion free as $Y \to \Sigma$ is a non-trivial bundle.

The trivial \spinc structure on $Y$ is the one with the trivial bundle of spinors $W^Y$. By viewing $W^Y$ as the
pull-back of a spinor bundle over $\Sigma$, one can endow it with a $\nabla^Y$-compatible spin connection that we
denote by $\nabla^Y$ as well.  The Clifford multiplication in $W^Y$ is as follows:  vectors in $H$ act
via the pull-back action, while $\vfi$ acts by $\pm i$ on $(W^Y)^\pm$, where the
splitting $W^Y=(W^Y)^+ \oplus (W^Y)^-$ is induced by the splitting of the spinor bundle over $\Sigma$.

Given a hermitian line bundle $E \to Y$ we say that the \spinc structure on $Y$ with the bundle of spinors $W_E^Y=W^Y
\otimes E$ is {\em determined } by $E$. A unitary connection $A$ in $E$ induces a spin connection $\nabla^Y
\otimes A$ in $W^Y \otimes E$; the  Dirac operator of this connection is denoted by $D_A$.  We need to understand
the moduli space of solutions of the perturbed 3-dimensional \swe\ on $Y$ (i.e., the equations defined using the above
Dirac operator) in a given \spinc structure. The space of reducible solutions in the \spinc structure determined by $E$
is non-empty only for torsion bundles $E$ (i.e., the ones with torsion Chern class). If $A_0$ is a smooth  flat
connection in $E$, then $A_0+\i \alpha$ is a reducible solution if and only if $\alpha$ is a closed one-form
on $Y$. This gives an identification  between the space of reducible solutions and the space of closed one-forms on
$Y$. The moduli space of reducible solutions, obtained by dividing the space of solutions by the gauge group action, is
therefore  identified  with $H^1(Y;\R)/H^1(Y)\iso H^1(\Sigma;S^1)$ via the choice of $A_0$.

Recall that the three-dimensional \swe\ are the equations for the critical points of the Chern-Simons-Dirac
functional on $\R \times Y$ (see \cite{mst}). This means that their linearization is self-adjoint, so for
positive dimensional moduli spaces the linearization of the equations is not surjective. The
appropriate notion of non-degeneracy of the moduli space is the following.
\begin{defn}
A component $\NN$ of the moduli space of \sws\ on $Y$ is (Morse-Bott) {\em non-degenerate}, if the kernel of
the linearization of the \swe\ at $(A,\Phi)$ is isomorphic to the tangent space to $\NN$ at $[A,\Phi]$ for any
point $[A,\Phi] \in \NN$.
\end{defn}

The following theorem, proved in \cite{moy},  describes the moduli spaces of solutions to the perturbed \swe\ on $Y$
for various \spinc structures.
\begin{theorem}\label{3-dim-sol}
Let $Y$ be a circle bundle of degree $n \ne 0$ over a surface $\Sigma$. The space of solutions in the \spinc structure
determined by $E \to Y$ is non-empty only if $E$ is the pull-back of a line bundle $F \to \Sigma$. Fix such an $E$ and let
$c= c_1(F)  \mod n$.\\
\odst{a}{The moduli space $\RR(E)$ of reducible solutions is homeomorphic to the dual torus $H^1(\Sigma; S^1)$.
Moreover,   if $c\not\equiv 0 \mod n$ or if $\Sigma\iso S^2$, this space is non-degenerate. In particular, the above
identification is  a diffeomorphism.}
\odst{b}{The irreducible components of the moduli space are parameterized by line bundles $E_0
\to \Sigma$ satisfying
$$c_1(E_0) \equiv c \mod n \hskip1.2cm \hbox{and}\hskip1.2cm 0<|c_1(E_0)| < g.$$
All of these components are compact and non-degenerate; they arise as the pull-backs of solutions  to the vortex
equations on $\Sigma$.}
\odst{c}{Since the Chern class of the \spinc structure is torsion,  the Chern-Simons-Dirac functional descends to a
real-valued function on the moduli space. If we normalize it so that it equals 0 on $\RR(E)$, then its value on the
component corresponding to a line bundle $E_0$ is $8\pi^2(c_1(E_0))^2/n$.}
\end{theorem}


\section{Seiberg-Witten solutions over a cylinder}
The first step in understanding the structure of the space of \sws\ on a manifold with a cylindrical end is
to study the solutions on the cylindrical part. A standard approach which guarantees good limiting behavior of
solutions at infinity is to consider only finite energy solutions (see \cite{mst} and \cite{mmr}).

Let $(A,\Psi)$ be a configuration on $[0,\infty) \times Y$ in a temporal gauge; denote by $(A_t,\Psi_t)$ the path
of configurations on $Y$ obtained by restricting  $(A,\Psi)$ to the slices $t \times Y$. Recall that the \swe\
on the cylinder take the form
$$\der{ }{t}(A_t,\Psi_t)=(*_Y(q(\Psi_t)-F_{A_t}),D_{A_t}\Psi_t).$$
The {\em energy}\ of a configuration $(A,\Psi)$ on $[0,\infty) \times Y$ in a tempo\-ral gauge is given by
the square of the $L^2$-norm of the right-hand side in the above equation.

Any solution on $Y$, being a critical point of the above equation, gives rise to a static solution on the cylinder;
such a solution clearly has finite energy. Moreover, any finite energy solution on the cylinder converges to a
static solution exponentially fast. This result is the Seiberg-Witten analogue of the exponential decay results
established by Morgan-Mrowka-Ruberman \cite{mmr} in Donaldson's theory. For our purposes, however, we do not need to
know that all solutions are exponentially decaying to solutions on $Y$. That is, without referring to exponential
decay results, we will consider only those configurations on the cylinder $[0,\infty) \times Y$ that decay
exponentially to solutions on $Y$, for some appropriately chosen decay constant.


\section{Seiberg-Witten equations over a cylindrical end manifold}
Let  $Z_0$ denote a compact oriented 4-manifold with boundary, $Y$, a circle bundle of degree $n \ne 0$  over a
surface $\Sigma$. Choose a collar $ [-1,0] \times Y \subset Z_0$ and equip $Z= Z_0 \cup [0,\infty) \times
Y$  with a cylindrical end metric that agrees with $dt ^2+g_Y$ on $[-1,\infty) \times Y$, where $g_Y$ is the
metric on $Y$ described in  Section \ref{circle-bundle}. As the  background connection $\nabla^Z$ for the Dirac
operator we use a metric compatible connection that agrees with the Levi-Civita connection on the complement of
$(-1,\infty) \times  Y$ and  agrees with the pull-back of the  connection $\nabla^Y$ on the cylinder $ [0,\infty)
\times  Y$.

Given a \spinc structure on $Z$ we denote the corresponding  bundles of spinors by $W=W^+ \oplus W^-$ and the
determinant line by $L=\det(W^+)$. As the configuration space for the \swe\ on $Z$ we choose the subset of
uniformly {\em exponentially decaying configurations}. The restriction of an exponentially decaying configuration
to the end $[0,\infty) \times Y$ differs from some static solution on the cylinder by a term that converges to zero
exponentially fast along the cylinder. We will specify the rate of convergence in Proposition \ref{Fredholm-property}.

Suppose now that we are in the situation from Section \ref{problem}: $Z$ is the complement of a representative
$\Sigma$ of a multiple class $d\xi$ in a closed manifold $X$. In this case the following proposition shows that any
\spinc structure on $Z$ which admits exponentially decaying solutions arises as the restriction of a \spinc structure
on $X$. We will use this to express the dimension of the moduli space of \sws\ on $Z$ in terms of the invariants of $X$.
\begin{proposition}\label{disk-extension}
Let $\Sigma \subset X$ be a smooth representative of the class $d\xi$ and let $n=(d\xi)^2$; denote by $N$ a compact
tubular neighborhood of $\Sigma$ and by $Z=X-\Sigma$ the cylindrical end manifold with the end modeled on $Y$.
Let $\SS$ be a \spinc structure on $Z$ for which the space of exponentially decaying \sws\ is non-empty.\\
\odst{a}{The restriction of $\SS$ to $Y$ is determined  by a torsion line bundle and it extends to a
\spinc structure on the disk bundle $N$. We fix the extension $W^N$ of the trivial \spinc structure
$W^Y$ on $Y$ to a \spinc structure on $N$, which is uniquely determined by requiring that the Chern class of its
determinant line is equal to $n[\Sigma]^*$. Through this choice we obtain a canonical extension of the given \spinc
structure $\SS$ on $Z$ to a \spinc structure on $X$.}
\odst{b}{Any two \spinc structures on $X$ that induce the given \spinc structure $\SS$ on $Z$ differ by a power of
the line bundle on $X$ with the Chern class $d\alpha$,  where $\alpha$ denotes the Poincar\'e dual of $\xi$.}
\end{proposition}
\begin{proof}
The induced \spinc structure on $Y$ is determined by a line bundle $E \to Y$ satisfying $L|_Y=E^2$. Note that the
existence of exponentially decaying \sws\ on $Z$ implies the existence of solutions on $Y$. Combining this with
Theorem \ref{3-dim-sol} we conclude that the line bundle $E$ is torsion and hence the pull-back of a line bundle on
$\Sigma$.

An extension of the trivial \spinc structure $W^Y$ on $Y$ to a \spinc structure on $N$ is given as follows. Denote
by $\hatN \to \Sigma$ the normal bundle of $\Sigma$ in $X$, considered as a complex line bundle of degree $n$. As a
complex manifold, the line bundle $\hatN$ carries a canonical \spinc structure with  the bundles of spinors
$W^+=\Lambda^0\oplus \Lambda^{0,2}$ and $W^-=\Lambda^{0,1}$. These bundles are determined up to isomorphism
by their restrictions to the zero-section $\Sigma \subset \hatN$. Note that $W^-|_\Sigma\iso \hatN \oplus
K_\Sigma^{-1}$, where $K_\Sigma$ is the canonical bundle of $\Sigma$. The  determinant line is therefore isomorphic
to $\hatN \otimes K_\Sigma^{-1}$. Since  the pull-back of $\hatN$ is trivial over $Y$, we need to change the \spinc
structure by a square root of the pull-back of the canonical bundle of $\Sigma$.

A different extension of a \spinc structure on $Z$ to a \spinc structure on $X$ can be obtained by changing the
\spinc structure on $N$ by a power of the pull-back of $\hatN$ (since this operation preserves the \spinc structure
on $Y$). From $c_1(\hatN)= (d\xi)^2[\Sigma]^*$ and Proposition \ref{cohomology-groups} it follows that the Chern class of
the auxiliary bundle on $X$ changes under this operation by a multiple of $d\alpha$. To see that these are the only
possibilities, consider two \spinc structures on $X$ that differ by a line bundle $E$. If they restrict to give the same
\spinc structure on $Z$, then $E|_Z$ is trivial, hence $c_1(E)$ lies in the kernel of the restriction homomorphism
$H^2(X) \to H^2(Z)$. Recall that this kernel is generated by $d\alpha$.
\end{proof}

As in the case of \spinc structures on $Y$ we will say that a \spinc structure on $N$ is {\em determined} by a line
bundle $E \to N$ (or by a line bundle $E_0  \to \Sigma$) if the bundle of spinors is of the form $W^N \otimes E$ (or
$W^N \otimes p^*E_0$) .

Given a \spinc structure on $X$, denote its determinant line by $Det$. With the above notation, we write
$\langle c_1(Det),\xi\rangle =k+(2s+1)d\xi^2$ for some $k \in \{0,\ldots,2d\xi^2-1\}$ and $s \in \Z$.
Note that $k$ is the representative of the residue class
$$\langle c_1(Det),\xi\rangle+d\xi^2 \mod 2d\xi^2$$
in $\{0,\ldots,2d\xi^2-1\}$; the purpose of the shift by $d\xi^2$ is to make $k$ directly related to the induced \spinc
structure on $Y$. Indeed, the induced bundle of spinors over $N$ is determined by a line bundle $E_0 \to \Sigma$ with
$c_1(E_0)=e$ satisfying $2e=d(k+2sd\xi^2)$.

Possible values of $k$ are constrained by the fact that $c_1(Det)$ is a characteristic class; this motivates the
following definition.
\begin{defn}\label{k-characterictic}
We call $k \in \{0,\ldots,2d\xi^2-1\}$ a {\em characteristic number} for $(X,\xi,d)$, if there exists a \spinc structure on
$X$ whose determinant line  $Det$ satisfies $\langle c_1(Det),\xi\rangle =k+(2s+1)d\xi^2$ for some $s \in \Z$.
\end{defn}

Clear the parity of a characteristic number is uniquely determined. It is easy to verify
that the parity is the only condition as described below.
\begin{lemma}\label{characteristic-conditions}
Let $X$ be a closed four-manifold, $\xi \in H_2(X)$ a primitive class of positive self-intersection, and $d\ge 1$ an
integer. If $\xi^2$ is even, then the set  of characteristic numbers for $(X,\xi,d)$ consists of all even numbers in
$\{0,\ldots,2d\xi^2-1\}$; if $\xi^2$ is odd, then the set of characteristic numbers consists of all numbers in
$\{0,\ldots,2d\xi^2-1\}$ with the  parity opposite to that of $d$. In particular, $kd$ is even in either case.
\end{lemma}


\section{Configuration space over a cylindrical end manifold}
In order to be able to consider the  \swe\ on $Z$ as an elliptic system of equations we need to choose
appropriate function spaces in which to study the equations. Since we chose to work with configurations on $Z$ that
uniformly exponentially decay to configurations on $Y$, we can
work  with the weighted Sobolev spaces $L^2_{r,\delta}$ for some small $\delta >0$. We cannot choose
$\delta=0$ because with this choice  the operator on $Y$, associated to the linearization of the \swe\ on  $Z$, has
a non-trivial  kernel and so the equations are not Fredholm (cf.~\cite{lockhart-mcowen}). Recall that the
$L_\delta^{2}$ norm of a function $f$ is defined by  $||f||_{\delta}^2=\int_Z |f|^2e^{\delta\tau}$, where $\tau \col Z
\to \R$ is a smooth function that is equal to $-1$ on the complement of $(-1,\infty) \times Y$ and agrees  with the $t$
coordinate on the  cylinder $[0,\infty) \times Y$; we further assume that $\tau$ depends only on $t$ and is
non-decreasing. Sobolev norms for $r>0$ are defined analogously.

We are interested in the space of configurations on $Z$ decaying to reducible solutions on $Y$. Fix for now a \spinc
structure on $Z$ with non-empty space of exponentially decaying \sws\ with reducible limits. For $(A,\Psi)$ a smooth
solution on $Z$ denote its asymptotic value by $(B_\infty,0)$. The space of all asymptotic values is an affine space
whose underlying vector space is the space of imaginary-valued closed one-forms on $Y$. We let $B$ be a smooth unitary
connection in $L$ (the determinant line of the \spinc structure on $Z$) which agrees with the pull-back of $B_\infty$
on the cylinder $[0,\infty) \times Y$. We use the connection $B$ (resp.~$B_\infty$) to identify the space of unitary
connections in $L$ (resp.~$L|_Y$) with the imaginary-valued 1-forms on $Z$ (resp.~$Y$).

Rather than working with the full configuration space, we restrict the possible asymptotic values of configurations
(by fixing the gauge at infinity). Specifically, we replace the space of imaginary-valued closed one-forms on $Y$ by
the subspace $\HH$ of imaginary-valued harmonic one-forms. Note that the subgroup of the gauge group $\GG^Y$ that
acts on $\HH$ consists of harmonic gauge transformations and is therefore isomorphic to $S^1 \times H^1(Y)$,
where $S^1$ denotes the constant gauge transformations. However, apart from the constant gauge transformations, the
only harmonic gauge transformations on $Y$  that extend to gauge transformations on $Z$ are those that correspond to
the classes in the image of $H^1(Z) \to H^1(Y)$. In our application of the cylindrical end moduli space to the genus
problem we will assume that $b_1(Z)=0$. Assuming the latter, the gauge group on $Z$ for the restricted configuration
space consists of the gauge transformations that converge to the constant gauge transformations on $Y$. We will
see later that the space of solutions on $Z$ contains reducible configurations, so the subgroup of constant gauge
transformations does not act freely on it. For this reason we consider the {\em based moduli space}, obtained by
dividing the space of solutions by the action of the gauge group based at infinity. The group of
constant gauge transformations still acts on the  based moduli space and we will use this action to obtain our results.
\begin{defn}\label{config-space}
The exponentially decaying {\em configuration space} on $Z$, corresponding to the reducible solutions on $Y$, is
defined to be $$\HH \times L^2_{2,\delta} (\i\Lambda^1(Z) \oplus W^+).$$ More precisely, the {\em configuration}
associated to an element $(h,\alpha,\Phi)$ is $$(A,\Psi)=(B+\alpha+\dtau h,\Phi).$$ A gauge transformation
$\sigma \in L_{3,loc}^2(Z,S^1)$ belongs to the {\em gauge group based at infinity}, $\GG_\infty$, if there exist
a $T>0$ and an $f \in L_{3,\delta}^2([T,\infty)\times Y,\i\R)$  so that the restriction of $\sigma$  to the cylinder
$[T,\infty) \times Y$ is given by $\sigma=\exp(f)$.
\end{defn}


\section{Deformation complex}
Recall that \spinc structures on $Z$ with non-empty exponentially decaying configuration space are induced
from $X$. Moreover, we will consider only \spinc structures on $X$ for which the induced \spinc structure on $Y$ is
non-trivial, unless $\Sigma$ is a sphere; this way the moduli space of reducible solutions on $Y$ is always
non-degenerate (cf.~Theorem \ref{3-dim-sol}). Suppose $(B_\infty,0)$ is the asymptotic value of a smooth solution
$(A,\Psi)$ on $Z$. We use the configuration space on $Z$ as described in Definition \ref{config-space}. Since the
configurations on $Z$ converge to reducible solutions on $Y$, \swe\ on $Z$ give rise to a map
\begin{equation}\label{sw-map}
\SW \colon \HH \oplus L^2_{2,\delta} (\i\Lambda^1(Z) \oplus W^+) \to L^2_{1,\delta}(\i\Lambda^{2,+}(Z) \oplus
W^-).
\end{equation}
This is well defined because all the terms in the \sw\ map exponentially decay to 0 (recall that the multiplication
$L^2_{2,\delta} \otimes L^2_{2,\delta} \to L^2_{1,\delta}$ is continuous). The {\em deformation complex}
$\DD_{(A,\Psi)}$ of the solution $(A,\Psi)$, taking into account the action of $\GG_\infty$, is
\begin{equation}\label{deformation-cx}
0\to L^2_{3,\delta}(Z,\i\R)  \nad{K_{(A,\Psi)}} \HH \oplus L^2_{2,\delta} (\i\Lambda^1(Z) \oplus W^+) \hbox to1.5cm{ }
\end{equation}
$$\hbox to3cm{ } \nad{T_{(A,\Psi)}\SW} L^2_{1,\delta}(\i\Lambda^{2,+}(Z) \oplus W^-)\to 0,$$
where $K_{(A,\Psi)}(f)=(0,2df,-f\Psi)$ is the infinitesimal gauge group action and
$$T_{(A,\Psi)}\SW (\alfahat,\psi)=(d^+\alfahat-2Q(\Psi,\psi),D_A\psi+\alfahat\cdot\Psi)$$
is the linearization of the \sw\ map at $(A,\Psi)$. Here $Q$ is the bilinear map associated to the quadratic map $q$ in
the \swe\, and  $\alfahat=\alpha+\dtau h$ is the one-form on $Z$ corresponding to $(h,\alpha) \in  \HH \oplus
L^2_{2,\delta} (\i\Lambda^1(Z))$. The cohomology groups of this complex provide some local information about the
based moduli space as described below. We first make the following observation.
\begin{lemma}
The zeroth cohomology group of the deformation complex is trivial.
\end{lemma}
\begin{proof}
If $f \in L^2_{3,\delta}(Z,\i\R)$ is in the kernel of $K_{(A,\Psi)}$, then $df=0$. Thus $f$ is constant and since it
converges to 0 at infinity, it must be identically equal to zero.
\end{proof}
The first cohomology group of the deformation complex (\ref{deformation-cx}) is called the {\em Zariski tangent
space} of the moduli space and the second cohomology group is called the {\em obstruction space}. If  $(A,\Psi)$ is a
regular point for the \sw\ map, then the obstruction space vanishes and the first cohomology of the complex is
isomorphic to the (geometric) tangent space of the moduli space at $[A,\Psi]$.

We will compute the index of the deformation complex $\DD_{(A,\Psi)}$ via the index of the {\em fiber complex}
$\FF_{(A,\Psi)}$ associated to it; the latter is defined via the following exact sequence of complexes
$$0 \to \FF_{(A,\Psi)} \to  \DD_{(A,\Psi)} \to \HH \to 0,$$
where $\HH$ denotes the deformation complex of the asymptotic value $(B_\infty,0)$ of $(A,\Psi)$ and the morphism
to $\HH$ corresponds to taking limits at infinity. Here we identified the complex $\HH$ with its only non-zero group
(in dimension 1), namely the group of harmonic one-forms on $Y$. The fiber complex differs from the full
deformation complex by a finite dimensional space $\HH$ (of dimension $2g$, where $g$ is the genus of $\Sigma$),
hence it suffices to compute the index of the fiber complex. The Fredholm properties of the fiber complex
are determined by the asymptotic behavior of its `wrapped-up' form (cf.~\cite{lockhart-mcowen}), given by
(\ref{fiber-cx}) below.
\begin{proposition}\label{Fredholm-property}
There exists $\delta_0>0$ so that for any $\delta \in (0,\delta_0]$, the \sw\ map (\ref{sw-map}), considered as a map
of  the spaces of $L^2_\delta$ sections, is Fredholm and its index is constant on the configuration space. Moreover,
if we conjugate the linearization to a map of the spaces of $L^2$ sections, then the kernel of its asymptotic map is
trivial.
\end{proposition}
\begin{proof}
Let $(A,\Psi)$ be a configuration on $Z$ with a reducible asymptotic value $(B_\infty,0)$ on $Y$. The $L^2_{\delta}$
adjoint $K_\delta^\ast$ of the infinitesimal gauge group action $K_{(A,\Psi)}$ is defined with respect to the
following inner products: for imaginary-valued forms $\alpha$ and $\beta$ let
$$\langle\alpha,\beta\rangle_\delta=\int_Z \alpha \wedge \ast \beta\, e^{\delta\tau},$$
where $\ast$ is the complex anti-linear extension of the Hodge star-operator;
for spinors $\psi$ and $\phi$ let
$$\langle\psi,\phi\rangle_\delta=2Re\int_Z \langle\psi(z),\phi(z)\rangle e^{\delta\tau} vol_Z.$$
Then $K_\delta^\ast(\alpha,\psi)=2e^{-\delta\tau}d^\ast e^{\delta\tau}\alpha+2\i Im\langle\Psi,\psi\rangle$; we
can drop the factor 2, thus obtaining the `wrapped-up' fiber complex
\begin{equation}\label{fiber-cx}
F_\delta \col L^2_{2,\delta}(\i\Lambda^1(Z) \oplus W^+) \to L^2_{1,\delta}(\i\Lambda^0(Z) \oplus \i\Lambda^{2,+}(Z)
 \oplus W^-),
\end{equation}
$$F_\delta(\alpha,\psi) = \big(e^{-\delta\tau}d^\ast e^{\delta\tau}\alpha+\i Im\langle\Psi,\psi\rangle,
d^+\alpha-2Q(\Psi,\psi), D_A\psi+\alpha\cdot\Psi\big).$$
To analyze this map we conjugate it by the isometry $T_\delta=e^{-\veps\tau} \col L^2 \to L^2_\delta$, where
$\veps=\delta/2$. This gives a map  $F$ between the spaces of $L^2$ sections that sends $(\alpha,\psi)$ to
$$\big(d^*\alpha-\veps\langle \alpha,d\tau\rangle+\i Im\langle\Psi,\psi\rangle,
d^+\alpha-\veps(d\tau\wedge \alpha)^+ -2Q(\Psi,\psi), D_A\psi-\veps d\tau\cdot\psi+\alpha\cdot\Psi).$$

For the purpose of computing the asymptotic operator of $F$ we only need to understand its form on the cylinder
$C=[0,\infty) \times Y$. Recall that $W^+ \iso p_2^*W_E^Y$ over $C$, where $W_E^Y$ is the bundle of spinors on $Y$
and $p_2 \colon [0,\infty) \times Y \to Y$ is the projection. Moreover, the Clifford multiplication by $dt$
induces an isomorphism between $W^+$ and $W^-$.  The bundles of forms on the cylinder are
given by $\Lambda^1(C)=p_2^*(\Lambda^0(Y) \oplus \Lambda^1(Y))$ and $\Lambda^{2,+}(C) \iso p_2^*\Lambda^1(Y)$, where
the last isomorphism follows from the fact that any self-dual two-form on the cylinder is of the form $dt \wedge
\gamma_t +\ast_3 \gamma_t$ for some path $\gamma_t$ of one-forms on $Y$. Writing $\alpha=\i fdt + \i\beta$,
$F(\i f,\i\beta,\psi)$ is given by
$$\Bigg(-\i\Big(\der{f}{t} -d_3^*\beta+\veps f\Big) ,
\i\Big(\der{\beta}{t} - d_3f+\ast_3d_3\beta- \veps\beta\Big), \der{\psi}{t} -D_{B_\infty}\psi-\veps\psi  \Bigg) + o(1).$$
Up to an obvious isomorphism, $F$ is of the form $\der{ }{t}-G+o(1)$, where the {\em asymptotic operator} $G$ acts on
the space of sections of $\Lambda^0(Y) \oplus \Lambda^1(Y) \oplus W_E^Y$ via the matrix operator
$$\left[\begin{matrix} -\veps & d_3^* & 0\cr  d_3 & -\ast_3d_3+\veps & 0\cr   0 & 0 &D_{B_\infty}+\veps\cr
\end{matrix}\right].$$
Notice that $G$ splits as the sum of (the perturbations  of) the asymptotic ope\-rators corresponding to the
anti-self-duality (ASD) operator and the Dirac operator on $Y$. By results of Lockhart and McOwen
\cite{lockhart-mcowen}, the operator $F$ (and hence $F_\delta$) is Fredholm if the kernel of $G$ is trivial. For the
ASD part this follows from the computation of the spectrum of this operator (see the proof of Proposition
\ref{fiber-index}): for $\veps=0$ the ASD asymptotic operator has non-trivial kernel, whereas for positive $\veps$
the kernel is trivial. For the Dirac part recall that by our choice of the \spinc structure on
$X$, the space of reducible solutions on $Y$ in the induced \spinc structure is non-degenerate, hence the kernel of
$D_{B_\infty}$ is trivial.  Moreover, the spectrum of $D_{B_\infty}$ depends only on the gauge equivalence class of
$B_\infty$; from compactness of the fundamental domain for the action of the gauge group on the space of flat
connections on $Y$, invertibility of the operators, and the fact that $D_{B_\infty}$ has discrete spectrum, it follows
that $G$ has trivial kernel for all small enough positive $\veps$.

Above we showed that the operator $F$ (and hence $F_\delta$), associated to an arbitrary configuration $(A,\Psi)$, is
Fredholm. To finish the proof we only need to show that the operator $F_\delta$ depends continuously on the
configuration $(A,\Psi)$. Suppose $(A',\Psi')$ is another configuration; denote the difference $(A',\Psi')-(A,\Psi)$
by $(a,\phi)$. Then the difference of the two linearizations is a bilinear map in $((a,\phi),(\alpha,\psi))$. The
required continuity now follows from the continuity of the  multiplication
$L^2_{2,\delta} \times L^2_{2,\delta} \to L^2_{1,\delta}$.
\end{proof}


\section{Index of the deformation complex}
We use the Atiyah-Patodi-Singer index formula (see \cite{aps:I}) to compute the index of the fiber complex.
This requires the operator to be independent of the $t$ variable along the cylinder, which is not the case at a
solution to the \swe. However, according to Lemma \ref{Fredholm-property}, the index can be computed using any
configuration with a reducible limit, in particular one which agrees with the pull-back of a reducible configuration
on $Y$ along the cylinder. Recall that the Atiyah-Patodi-Singer index can be expressed in terms of the (extended)
$L^2$ solutions;  since by Lemma \ref{Fredholm-property} the asymptotic operator on $L^2$ has no kernel, the
Atiyah-Patodi-Singer index of the fiber complex agrees with its Fredholm index.

We first recall the relevant results from \cite{aps:I}.  Let $E_0,\, E_1 \to Z_1$ be hermitian vector bundles of
the same fiber dimension over $Z_1=\{z \in Z \, | \, \tau(z) \le 1\}$ and let $F \col \Gamma(E_0) \to \Gamma(E_1)$
be an elliptic operator. Assume that on the cylinder $[0,1] \times Y \subset Z_1$ we have  $F=\sigma(\der{}{t}-G)$,
for some self-adjoint operator $G \col \Gamma(E) \to \Gamma(E)$, where $E=E_0|_Y$ and $\sigma$ is an isomorphism
between $E_0$ and $E_1$ over the cylinder. The new ingredient in the index formula is a boundary correction term,
which is a spectral function of $G$. More precisely, let $\eta(s)=\sum\limits_{\lambda \ne 0}
\sign(\lambda)|\lambda|^{-s}$, where $\lambda$ runs over the spectrum of $G$, be the {\em eta function} of $G$. This
series defines a holomorphic function in a half-plane $Re(z)>z_0$ and extends to a meromorphic function on the
entire plane; this extension has a finite value at 0. The correction term is defined in terms of
$\eta(0)$ and the dimension $h$ of the kernel of $G$.

The Atiyah-Patodi-Singer domain of $F$, denoted by $\Gamma(E_0,P)$, consists of all the sections of $E_0$ whose
restriction to the boundary $1 \times Y$ of $Z_1$ lies in the kernel of $P$. Here $P=P_{\ge 0} \col \Gamma(E) \to
\Gamma(E)$ denotes the spectral projection of $G$, corresponding to the non-negative eigenvalues, i.e., the orthogonal
projection onto the subspace spanned by the non-negative eigenvalues of $G$. Notice that with this choice of the
domain any solution of $F=0$ extends to an $L^2$ solution on $Z$.
\begin{theorem}\label{APS-index}
With the above notation, $F \col \Gamma(E_0,P) \to \Gamma(E_1)$ has a finite index given by
$$\ind_{APS}F=\int_{Z_1} k - \frac{h+\eta(0)}{2},$$
where $k$ is a differential form on $Z_1$ determined by $F$, called the {\em index density} of $F$. Moreover,
$\ind_{APS}F=h(E_0)-h(E_1)-h_\infty(E_1)$, where $h(E_0)$ is the dimension of the kernel of $F$ on the space
of $L^2$ sections on $Z$, $h(E_1)$ is the corresponding dimension for the adjoint $F^*$ of $F$, and $h_\infty(E_1)$ is
the dimension of the space of asymptotic values of extended $L^2$ solutions of $F^*$.
\end{theorem}

To compute the index of the fiber complex we need the signature eta invariant of $Y$, which was
computed by Komuro \cite{komuro}, and the eta invariant of the perturbed Dirac operator on $Y$ (see Section
\ref{circle-bundle}), which was computed by Nicolaescu \cite{nicolaescu:eta}; the latter
paper also contains a derivation of the index formula presented below.
\begin{proposition}\label{fiber-index}
Let $X$ be a closed four-manifold and $\Sigma \subset X$ a smooth surface representing the homology class $d\xi$,
where $d>1$ and $\xi \in H_2(X)$ is a primitive class of positive self-intersection. Given a \spinc structure on
$X$ with the determinant line $Det$, we write $\langle c_1(Det),\xi\rangle =k+(2s+1)d\xi^2$ for some characteristic
number $k$ and some $s \in \Z$. The index of the fiber complex (\ref{fiber-cx}), associated to the space of \sws\ on
the cylindrical end manifold $Z=X-\Sigma$, that along the end converge to a fixed reducible solution $(B_\infty,0)$, is
\begin{equation}\label{fiber-ind}
 \frac{1}{4}\int_Z c_1(B)^2 -\frac{\sigma(X)}{4} +\frac{(k-d\xi^2)^2}{4\xi^2}+1+b_1(X)-b_2^+(X)-2g  ,
\end{equation}
where $B$ is a unitary connection in $Det|_Z$, which agrees with the pull-back of $B_\infty$ on the end $[0,\infty)
\times Y$.
\end{proposition}
\begin{remark}
By rewriting the above dimension formula, one can obtain the Fr{\o}yshov invariant (see \cite{froyshov}) of the
circle bundle $Y$ for a range of \spinc structures.
\end{remark}
\begin{proof}
We will compute the (real) index of the linearization of the \sw\ map at a configuration $(B,0)$ with asymptotic
value $(B_\infty,0)$. In this case the associated operator $F$ (defined in the proof of Proposition
\ref{Fredholm-property}) on the spaces of $L^2$ sections takes the form
$$F(\alpha,\psi) \mapsto \big(d^*\alpha,d^+\alpha, D_B\psi\big)
-\veps\big(\langle \alpha,d\tau\rangle,(d\tau\wedge \alpha)^+, d\tau\cdot\psi \big).$$
Clearly $F$ splits as the sum $F_0 \oplus F_1$, where $F_0$ is a zeroth-order perturbation of the anti-self-duality
operator $\AA=d^*\oplus d^+$,  and $F_1$ is a zeroth-order perturbation of the Dirac operator $D_B$. Hence we can
split the index computation accordingly.

{\em Index of the anti-self-dual part:}
For the purpose of invoking the Atiyah-Patodi-Singer index theorem we complexify the spaces of forms. The index
density in the statement of
Theorem \ref{APS-index} depends only on the principal symbol of the operator; for $F_0$ it is therefore determined
by $\AA$. The difference in the indices of $F_0$ and $\AA$ comes from the correction term  in the index formula and
can be described as the spectral flow of a family of associated asymptotic operators as made precise below. Notice
that $\AA$ is (isomorphic to) the adjoint of the operator $A_+$ from \cite{aps:I}  and we will  use the following fact
from there. The Atiyah-Patodi-Singer index computation gives
$$\frac{1}{2}(\sigma(Z)+\chi(Z))=\int_Z k - \frac{1}{2}\eta_{\sign}(0),$$
where $\sigma(Z)$ and $\chi(Z)$ are the signature and the Euler characteristic of $Z$ respectively, $k$ is the index
density of $\AA$ and $\eta_{\sign}(0)$ is its eta invariant. Therefore we obtain
$$\ind_{APS}(\AA)=-\int_Z k -\frac{1}{2}(h-\eta_{\sign}(0))=-\frac{1}{2}(\sigma(Z)+\chi(Z)+h),$$
where $h$ is the dimension of the kernel of the asymptotic operator
$$G_0=\left[\begin{matrix} 0 & d_3^* \cr  d_3 & -\ast_3d_3 \end{matrix}\right] \col \Omega^0(Y) \oplus
\Omega^1(Y) \to \Omega^0(Y) \oplus \Omega^1(Y).$$
A pair $(f,\beta)$ in the kernel of $G_0$ satisfies $d_3^*\beta=0$ and $ d_3f = \ast_3d_3\beta=0$. Hence the kernel
of $G_0$ consists of harmonic forms and we have $h=1+2g$. Using this, along with $\sigma(X)=\sigma(Z)+1$,
$\chi(X)=\chi(Z)+2-2g$, and $\sigma(X)+\chi(X)=2-2b_1(X)+2b_2^+(X)$, we obtain $\sigma(Z)+\chi(Z)=
-1+2b_1(X)+2b_2^+(X)+2g$; thus $\ind_{APS}(\AA)= b_1(X)-b_2^+(X)-2g$.

Notice that the asymptotic operator of $F_0$ is of the form $G_0+\veps E$, where $E=\left[\begin{matrix} -1 & 0
\cr  0 & 1 \end{matrix}\right]$. It will be convenient to consider the family of operators $G_u=G_0+u\veps E$ for
$u \in [0,1]$, connecting the limiting operators $G_0$ and $G_1$ of $\AA$ and $F_0$ respectively. The difference in
the correction terms for $G_0$ and $G_1$ is equal to the spectral flow of the family $G_u$. More precisely, let
$\etahat_u=\frac{1}{2}(h_{G_u}+\eta_{G_u}(0))$ denote the {\em reduced eta invariant} of $G_u$. From the definition
it is  clear that $\etahat_u$ has a jump at $u=a$ only if some eigenvalue $\lambda(u)$ of $G_u$ vanishes at $a$ (but
is not zero at least on one side of $u=a$). The case of interest to us is when any eigenvalue $\lambda(u)$ behaves in
one of the following ways: $\lambda(u)=0$ for all $u$, $\lambda(u)\ne 0$ for all $u$, or $\lambda(u)=0$ iff $u=0$ and
this zero is transverse.  Assume first that $\lambda(u)$ is the only eigenvalue crossing 0 at $u=0$. If
$\lambda'(0)>0$, then $\etahat_1=\etahat_0$, since at $u=0$ the eigenvalue $\lambda(0)$ contributes $+1$ to $h$ and
nothing to $\eta$, whereas for $u>0$ it contributes $+1$ to $\eta$ and nothing to $h$.  Similar considerations in the
case $\lambda'(0)<0$ imply that $\etahat_1=\etahat_0-1$, since in this case $\lambda(u)$ contributes to $h$ and
$\eta$ with the opposite signs. This clearly generalizes to a finite number of eigenvalues crossing 0 at $u=0$.

To determine the difference of indices of $\AA$ and $F_0$ we therefore need to understand the behavior of the
eigenvalues of the family $G_u$. A pair $(f,\beta) \in \Omega^0(Y) \oplus \Omega^1(Y)$ in the kernel of $G_u$
satisfies
$$d_3^*\beta=\veps u f, \ \ d_3f-\ast_3d_3\beta=-\veps u \beta.$$
This implies $d_3^*d_3f=-u^2\veps^2f$, which has no solutions for $u \ne 0$ since the Laplace operator $d_3^*d_3$
is positive definite. For $u=0$ we computed above that the dimension of the kernel is $1+2g$. The eigenvalue
$\lambda_1(u)=-\veps u$, corresponding to the space of constant functions,  has multiplicity 1 (at 0), whereas the
multiplicity of $\lambda_2(u)=\veps u$, corresponding to the space of harmonic one-forms,  is $2g$. We conclude
from the previous paragraph that $\etahat_1=\etahat_0-1$, and hence the index of $F_0$ is $1+b_1(X)-b_2^+(X)-2g$.

{\em Index of the Dirac part:}
Using similar considerations as above we see that the index of $F_1$ is equal to the index of $D_B$;  the reason for
this is that all the asymptotic operators $D_{B_\infty}+u\veps$ are invertible (for $u \in [0,1]$) by Theorem
\ref{3-dim-sol} and our choice of $\veps$. The index density of $D_B$ is $\frac{1}{8}(c_1(B)^2-L(\nabla^Z))$, where
$L(\nabla^Z)=\frac{1}{3}p_1(\nabla^Z)$ is the Hirzebruch $L$-class associated to the  background connection
$\nabla^Z$ in $Z$; recall that $\nabla^Z$ agrees with the Levi-Civita connection $\nabla^{LC}$ on the complement of
$(-1,\infty) \times Y$ and agrees with the pull-back of $\nabla^Y$ on $[0,\infty) \times Y$. From Theorem
\ref{APS-index} we get $\ind_{APS}(D_B)=\frac{1}{8}\int_Z(c_1(B)^2-L(\nabla^Z))-\frac{1}{2}
\eta_{D_{B_\infty}}(0)$. The eta invariant of the Dirac operator on $Y$, coupled to the flat connection
$B_\infty$, was computed in \cite{nicolaescu:eta} and is equal to
$$ \eta_{D_{B_\infty}}(0)=-\frac{ad^2}{6}-\frac{k^2}{4a}+\frac{kd}{2},$$
where we used the fact that the radius of the fiber circles in $Y$ is 1. The only other term we need to
interpret is the integral of the $L$-class. If we were using the Levi-Civita connection as the background
connection on $Z$, then we could use the fact that $L(\nabla^{LC})$ is the index density of the signature operator
on $Z$. In particular, we have $\sigma(Z)=\int_Z L(\nabla^{LC}) -\eta_{\sign}(0)$. The signature eta invariant
$\eta_{\sign}(0)$ for $Y$ was computed in \cite{komuro} and is given by
$$\eta_{\sign}(0)=1-\frac{ad^2}{3}+\frac{2ad^2}{3}\big( a^2d^4+2g-2\big)  .$$
In our case, however, there is another term coming from the difference in the $L$-classes of the two connections:
we have
$$\int_{Z_0} L(\nabla^Z)  = \int_{Z_{-1}}  L(\nabla^{LC}) + \int_{[-1,0] \times Y} L(\nabla^Z).$$
The last term in this expression can be computed explicitly; write $\nabla^Z=\nabla^\infty+f(t)\alpha$, where
$\nabla^\infty$ denotes the pull-back connection, $\alpha=\nabla^\infty-\nabla^{LC}$ and $f$ is a smooth
non-decreasing function that maps $[-1,0]$ onto itself.  The computation can be done with respect to a local
orthonormal coframe $(\vfi_1,\vfi_2,\vfi)$ on $Y$, where $\i\vfi$ is the connection  of the circle bundle $Y \to
\Sigma$ and $(\vfi_1,\vfi_2)$ is the pull-back of a local coframe on $\Sigma$. This yields
$$\int_{[-1,0] \times Y} L(\nabla^Z)=-\frac{2ad^2}{3}\big( a^2d^4+2g-2\big)$$
(see \cite{nicolaescu:eta} for more details). Note that the index formula gives the complex index  of the operator
$F_1$.
\end{proof}


\section{Dimension of the cylindrical end moduli space}
To express the formal dimension of the moduli space of \sws\ on the cylindrical end manifold $Z=X-\Sigma$ in
terms of the data on the closed manifold $X$ we need the following result.
\begin{lemma}\label{c1^2-flat}
Let $X$ be a  closed four-manifold and let $\Sigma \subset X$ be an embedded surface with self-intersection $n \ne 0$.
Denote by $Z$ the complement $X-\Sigma$, thought of as a manifold with a cylindrical end $[0,\infty) \times Y$. Given
a \spinc structure on $X$, let $p=\langle c_1(Det),[\Sigma]\rangle \in \Z$, where $Det \to X$ denotes the determinant
line of the \spinc structure. For any  unitary connection $B$ on $L = Det|_Z$, whose restriction to the cylinder
$[0,\infty) \times Y$ agrees with the pull-back of a flat connection in $Det|_Y$, we have
\begin{equation}\label{c1^2-formula}
 \int_Z c_1(B)^2 = c_1(Det)^2 - \frac{p^2}{n}.
\end{equation}
\end{lemma}
\begin{proof}
Let $Z_1=\{z \in Z\, |\, \tau(z) \le 1\}$ and let $Y$ be the oriented boundary of $Z_1$. We think of $X$ as the union
of $Z_1$ and  a compact tubular neighborhood $N$ of $\Sigma$ in $X$. Denote by $[0,1] \times Y$ the oriented
collar to the boundary $Y$ of $N$.  Suppose $A$ is a connection in $Det|_N$ that in a neighborhood of the boundary
$0 \times Y$ agrees with
the pull-back of $B_\infty$, the latter being the limit of $B$.  Then $A$ and $B$ together define a (smooth)
connection in $Det$ and we have $c_1(Det)^2=\int_{Z} c_1(B)^2+ \int_N c_1(A)^2$. We will evaluate the second integral.
Combining Lemma \ref{Fredholm-property} and the index formula (\ref{fiber-ind}), we see that, for the purpose of the
computation, we can choose $A$ to be any connection in $Det|_N$ which in a neighborhood of $Y$ agrees with the
pull-back of some flat connection in $Det|_Y$. In what follows we  use the notation from Section \ref{circle-bundle}.

Let $A_1$ be the pull-back of a constant curvature connection in $Det|_\Sigma$ to $Det|_N$; then $F_{A_1}=
-2\pi\i p\omega$. We let $A=A_1+\i\frac{p}{n}f(t)\vfi$, where $i\vfi$ is a constant curvature connection
of the circle bundle $Y \to \Sigma$ and $f \col [0,1] \to [0,1]$ is a smooth non-increasing function which is
identically 1 in a neighborhood of 0 and identically 0 in a neighborhood of 1. Since the degree of $Y$ is $-n$, we have
 $d\vfi=2\pi n\omega$, so $A$ is flat close to the boundary of $N$. A simple computation
shows that $F_A=2\pi\i p(f(t)-1)\omega + \i \frac{p}{n}f'(t) dt\wedge \vfi$ and hence
$$ F_A \wedge F_A = -4\pi \frac{p^2}{n}\big(f(t)-1\big)f'(t)\, dt \wedge vol_Y.$$
Using this along with $vol(Y)=2\pi$ and $c_1(A)=\frac{\i}{2\pi}F_A$ gives
$$\int_Z c_1(A)^2 = 2\frac{p^2}{n}\int_0^1(f(t)-1)f'(t)dt=\frac{p^2}{n}.$$
\end{proof}

Now we can obtain a convenient formula for the (formal) dimension of the moduli space of \sws\ on $Z$ with reducible
limits.
\begin{corollary}\label{formal-dim}
Let $X$ be a closed four-manifold with $b_1(X)=0$ and $b_2^+(X)\allowbreak =1$. Suppose $\Sigma \subset X$ is a smooth surface
representing the homology class $d\xi$, where $d>1$ and $\xi \in H_2(X)$ is a primitive class of positive
self-intersection. Let $Z=X-\Sigma$, thought of as a cylindrical end manifold. Given a \spinc structure on
$X$, let $p=\langle c_1(Det),\xi\rangle$, where $Det \to X$ denotes the determinant line of the \spinc structure; we
write $p=k+(2s+1)d\xi^2$ for some characteristic number $k$ and some $s \in \Z$. The formal dimension of the based
moduli space of \sws\ on $Z$ with reducible limits is given by
\begin{equation}\label{formal-dimension}
\frac{c_1(Det)^2-\sigma(X)}{4} +\frac{(k-d\xi^2)^2-p^2}{4\xi^2}=\frac{c^2-\sigma(X)}{4}
+\frac{(k-d\xi^2)^2}{4\xi^2},
\end{equation}
where $c \in H^2(X;\R)$ is a class of non-positive square.
Moreover, this dimension depends only on the induced  \spinc structure on $Z$, i.e., it is independent of $s$.
\end{corollary}
\begin{proof}
Starting with the index of the fiber complex, given by (\ref{fiber-ind}), recall from the discussion preceding
Lemma \ref{Fredholm-property} that we need to add $2g$ to it (where $g$ is the genus of $\Sigma$) to obtain the index
of the full deformation complex.
The left-hand formula now follows from the lemma above with $n=(d\xi)^2$ and $p$ replaced by $pd$.

Considering $c_1(Det)$ as a class in $H^2(X;\R)$,  we write
$$c_1(Det)=\frac{p}{a}\alpha+c$$
for some $c\in H^2(X;\R)$, where $\alpha$ is the Poincar\'e dual of $\xi$ and $a=\xi^2$. It follows that $c \cup
\alpha= \langle c,\xi\rangle=0$, hence $c_1(Det)^2=p^2/a+c^2$, which implies the right-hand formula. Since $\alpha$
and $c$ are orthogonal and $b_2^+(X)=1$, $c^2$ cannot be positive.

Finally,  consider another \spinc structure on $X$ which induces the same \spinc structure on $Z$. Then its
determinant line, denoted by $Det'$, satisfies  $c_1(Det')=c_1(Det)+2sd\alpha$ for some $s \in \Z$ (cf. Proposition
\ref{disk-extension}). This shows that $c'$, defined analogously as $c$ above, equals $c$; the last assertion then
follows from the right-hand formula.
\end{proof}


\section{Compactness and regularity of the cylindrical end moduli space}\label{comp-reg}
Let $\Sigma$ be an embedded surface of positive self-intersection in a closed four-manifold $X$ with $b_1(X)=0$.
Denote by $Z$ the complement $X-\Sigma$, thought of as a manifold with a cylindrical end $[0,\infty) \times Y$.
We topologize the moduli space of exponentially decaying \sws\ on $Z$ by the weakest topology, containing the
topology of uniform $C^k$ convergence on compact subsets for some large $k$, with respect to which the
Chern-Simons-Dirac functional along the cylinder is continuous at infinity.
\begin{proposition}\label{compactness}
Fix a \spinc structure on $X$ for which there are no irreducible \sws\ on $Y$ in the induced \spinc structure.
Then  the moduli space of exponentially decaying \sws\ on $Z$ (in the induced \spinc structure) with reducible limits
is compact.
\end{proposition}
\begin{proof}
Given a sequence $(A_n,\Psi_n)$ of \sws\ on $Z$, observe that this sequence has a convergent subsequence on any
submanifold $Z_t \subset Z$ (after appropriate changes of gauge). This is essentially a consequence of
\cite[Lemma 4]{kronheimer-mrowka:thom}. The uniform boundedness of $|\Psi_n|$ on $Z$ follows from the fact that the
configurations converge to zero at infinity, and from a standard maximum principle argument
\cite[Lemma 2]{kronheimer-mrowka:thom}. The only difference is that we are not using the Levi-Civita connection as
the background connection. However, the two Dirac operators differ by Clifford multiplication by a
one-form $f(t)dt$, where $f$ is a smooth bounded function (cf.~\cite[Lemma 5.2.1]{moy}), and it easily
follows from this that the result holds for the perturbed Dirac operator as well.

Using the diagonal argument we can thus find a subsequence of $(A_n,\Psi_n)$ which, after appropriate changes of
gauge, converges on all compact subsets of $Z$ to some solution $(A,\Psi)$ of the \swe\ on $Z$; we still denote this
subsequence by $(A_n,\Psi_n)$. Potential non-compactness therefore arises from the behavior of solutions
on the end $[0,\infty) \times Y$. In particular, the convergence of the sequence $(A_n,\Psi_n)$
depends on the convergence of the sequence of its asymptotic values. Recall that the moduli space of asymptotic
values is identified with the space $\HH$ of imaginary valued harmonic one-forms on $Y$.

We first prove that the sequence of the asymptotic values of $(A_n,\Psi_n)$ is bounded; in fact, the moduli space is
contained in the fiber of the projection to the space of asymptotic values. Let $(A,\Psi)$ and $(B,\Phi)$ be two
solutions with reducible asymptotic values. Then $\alpha=A-B$ exponentially decays to a form $\alpha_\infty \in \HH$.
To prove that $\alpha_\infty=0$ it suffices to show that $\int_K \alpha_\infty=0$ for any embedded circle $K \subset Y$
representing a homology class in the kernel of the morphism $H_1(Y) \to H_1(Z)$. Since such a $K$ is the boundary of a
surface $S \subset Z$, this is equivalent to $\int_S (F_A-F_B)=0$. As both $F_A$ and $F_B$ represent the same
relative cohomology class on $Z$, the claim follows.

Suppose now that $[(A_n,\Psi_n)]$ does not converge to $[(A,\Psi)]$ in the topology of the moduli space. This means
that the Chern-Simons-Dirac functional has different values on the asymptotic configurations of $(A_n,\Psi_n)$ and
$(A,\Psi)$. Then the asymptotic value of $(A,\Psi)$ is irreducible (see \cite[Proposition
8.5]{mst} for more details), which contradicts the assumption that the space of solutions on $Y$
consists entirely of reducibles.
\end{proof}

Now we turn to the question of regularity of the moduli space. We are particularly interested in behavior of the Dirac
operator at a reducible solution; a more general result regarding such Dirac operators on a closed manifold is proved
in \cite{salamon}. For technical reasons we choose to work with $L_{4,\delta}^2$ configuration space.
\begin{proposition}\label{smoothness}
Let $X$ be a closed oriented four-manifold with $b_1(X)=0$ and $b_2^+(X)=1$. Suppose that $\Sigma \subset X$ is a
smooth surface of positive self-intersection, and let $Z=X-\Sigma$, thought of as a cylindrical end manifold.
Then for any small enough $\delta>0$ there exists a second category subset $\Omega$ of imaginary-valued
$L^2_{3,\delta}$ self-dual two-forms on $Z$, such that for any $\omega\in\Omega$ the following holds. For any
exponentially decaying connection $A$ in the determinant line $L \to Z$ satisfying $F_A^+=\omega$, the Dirac operator
$D_A$ is either injective or surjective. Moreover, the irreducible part of the $\omega$-perturbed  moduli space of
exponentially decaying \sws\ on $Z$ with reducible asymptotic values is a smooth orientable finite dimensional
manifold.
\end{proposition}
\begin{proof}
Fix a smooth unitary connection $A_0$ in the determinant line $L \to Z$ and let $\NN$ be the manifold of all
exponentially decaying configurations $(A,\Psi)$ (with reducible limits) satisfying
$$T(A,\Psi)=(d^*(A-A_0),D_A\Psi)=0, \ \ \Psi\ne 0.$$
To prove that $\NN$ is a smooth manifold we need to verify that the differential
$$D_{(A,\Psi)}T(\alpha,\psi)=(d^*\alpha,D_A\psi+\alpha\cdot\Psi)$$
of $T$ at $(A,\Psi)$ is onto. The adjoint of this operator is
$$(f,\chi)\mapsto(df+\i\langle\_ \cdot\Psi,\chi\rangle-\delta f\, d\tau ,D_A\chi+\delta d\tau\cdot\chi),$$
where we used the usual $L^2_{\delta}$ inner product on the space of imaginary valued one-forms and the real
part of the hermitian $L^2_{\delta}$ inner product on spinors. Expression $\i\langle\_ \cdot\Psi,\chi\rangle$
denotes the imaginary valued one-form characterized by
$$\langle\alpha,\i\langle\_ \cdot\Psi,\chi\rangle\rangle=Re\langle\alpha \cdot\Psi,\chi\rangle,$$
for all imaginary valued one-forms $\alpha$.

Rather than proving directly that the kernel of the adjoint is trivial, we will replace the adjoint by the operator of
the same kind with $\delta=0$. The injectivity of thus obtained operator $D$ implies the injectivity of the adjoint
for $\delta>0$ small enough.

Suppose that $D(f,\chi)=0$. Computing with respect to a local orthonormal coframe $\{\vfi_1,\ldots,
\vfi_4\}$ on $Z$ we obtain
$$\i d^*df=Re\langle\i D_A\Psi,\chi\rangle +\langle\i\Psi,D_A^*\chi\rangle
-Re\langle\i \sum_j (\nabla_{\vfi_j}\vfi_j+(*d*\vfi_j)\vfi_j)\cdot\Psi,\chi\rangle.$$
The first two terms in the last expression vanish as $\Psi$ and $\chi$ are harmonic. Moreover, the one-form
inside the last term vanishes for Levi-Civita connection. The background connection we are using differs
from the Levi-Civita connection on the cylinder $[0,\infty) \times Y$ by a multiple of
$$\left(\begin{array}{cccc}
0 & \vfi_3 & \vfi_2 & 0\cr
-\vfi_3 & 0 & -\vfi_1 & 0\cr
-\vfi_2 &\vfi_1 & 0 & 0\cr
0 & 0 & 0 & 0 \end{array}\right)$$
(see \cite[Lemma 5.2.1]{moy}), where we took $(\vfi_1,\vfi_2)$ to be the pull-back of an orthonormal coframe
on $\Sigma$ and $\i\vfi_3$
to be the connection one-form of the circle bundle $Y$. So the last term in the above expression vanishes
as well. This implies that $f$ is constant, and since $f$ exponentially decays to zero, $f=0$. Finally, since
$\langle\alpha\cdot\Psi , \chi\rangle=0$ for any $\alpha$, where $\Psi$ and $\chi$ are harmonic, it follows from
unique continuation property for harmonic spinors that $\chi=0$ (see \cite[Lemma 6.2.1]{morgan:swbook}). This proves
that $\NN$ is smooth. In particular, at any point $(A,\Psi) \in \NN$,
$$\im D_A + \{\alpha \cdot \Psi \mid d^*\alpha=0\}=L^2_{3,\delta}(W^-).$$

Let $\Omega_0$ be the set of regular values of the map $\NN \to L^2_{3,\delta}(\i\Lambda^{2,+}(Z))$, $(A,\Psi) \mapsto
F_A^+$. Given $\omega \in \Omega_0$, let $(A,\Psi) \in \NN$ be a point satisfying $F_A^+=\omega$. Now the differential
$d^+$ from the tangent space to $\NN$ at $(A,\Psi)$ to the space of imaginary self-dual two-forms is onto by choice
of $\omega$.  Further, from Hodge decomposition of $L^2_{k,\delta}$ forms on $Z$ (see \cite{lockhart:hodge}; note that
we chose $\delta$ so that the Laplace operator is Fredholm), and the assumption $H^1(Z)=0$, it follows that the space
of co-closed one-forms maps isomorphically onto the space of self-dual two-forms. Thus for any co-closed one-form
$\alpha$ there exists a spinor $\psi$ so that $D_A\psi+\alpha\cdot\Psi=0$.
Combining this with the observation at the end of the previous paragraph proves that $D_A$ is onto.

The proof of regularity of the irreducible part proceeds as for a closed manifold (see \cite{morgan:swbook}); we get
a set $\Omega_1$ of regular perturbation parameters and let $\Omega=\Omega_0 \cap \Omega_1$.
\end{proof}



\vskip2cm
\centerline{\large \sc Part II: Genus bounds}
\section{An example: genus bounds in $\cp^2$}
Before stating and proving the main theorem (in the next section), we will demonstrate the argument in the simplest
possible case, for $X=\cp^2$. Let  $\xi=[\cp^1]$ be the standard generator of $H_2(X)$; then $\xi^2=1$. We
fix $d>1$ and consider a smooth genus $g$ representative $\Sigma$ of the class $d\xi$. A \spinc structure on $X$ is
uniquely determined by $p=\langle c_1(Det),\xi\rangle$. As before we write $p=k+d$ for some characteristic
number $k$. We will see in the proof of Theorem \ref{main-theorem} that it suffices to
consider $k \in \{0,\ldots,d\}$. Recall from the discussion preceding Definition \ref{k-characterictic} that
the line bundle, determining the induced \spinc structure on $Y$, is the pull-back of a line bundle $E_0 \to \Sigma$
with $c_1(E_0)={kd/2}$. According to part (b) of Theorem \ref{3-dim-sol}, irreducible solutions in the given
\spinc structure on $Y$ exist only if
\begin{equation}\label{cp^2-condition}
g > \frac{kd}{2}.
\end{equation}

The general  formula (\ref{formal-dimension}) for the formal dimension of the based moduli space of  \sws\ with
reducible limits on $Z$  in the case we are considering becomes
\begin{equation}\label{cp^2-dimension}
\frac{(k-d)^2-1}{4}.
\end{equation}
Notice that this number is even and by 1 greater than the expected dimension of the moduli space.

Using inequality (\ref{cp^2-condition}) and the dimension formula (\ref{cp^2-dimension}) we obtain a lower bound
on the genus $g$ of $\Sigma$ based on the following observation (explained in the proof of Theorem
\ref{main-theorem}): if the moduli space of \sws\ with reducible limits on $Z$ is compact and positive dimensional,
this leads to a contradiction. In other words, using Proposition \ref{compactness}, if the moduli space is positive
dimensional, then (\ref{cp^2-condition}) must hold.

The dimension (\ref{cp^2-dimension}) of the based moduli space is positive for $k<d-1$ in the range of $k$'s
considered. Since $d$ and $k$ have different parity,  we
conclude (for $d \ge  3$) that
$$ g > \frac{(d-3)d}{2}. $$
Since the classes $d[\cp^1]$ for $d=1,\, 2$ are represented by spheres, this inequality is equivalent to the {\em
Thom conjecture}, which was first established by Kronheimer and Mrowka \cite{kronheimer-mrowka:thom}.
\begin{theorem}
Let $\Sigma \subset \cp^2$ represent the class $d[\cp^1]$ for some $d \ge 1$. Then the genus $g$ of $\Sigma$
satisfies
\begin{equation}\label{cp^2-bound}
g \ge \frac{(d-1)(d-2)}{2}.
\end{equation}
Moreover, this lower bound is attained by a smooth holomorphic curve representing this homology class.
\end{theorem}
\begin{remark}
(a) Note that the same genus bound holds for classes $d[\cp^1]$ with $d\le-1$, where $d$ gets replaced by $|d|$ in
(\ref{cp^2-bound}). This observation is true in general, since we can always replace $\xi$ by $-\xi$.\\
(b) From above computations we see that genus bound (\ref{cp^2-bound}) holds in any $X$ which is a
rational homology $\cp^2$; by possibly changing the orientation we can assume that $X$ is positive definite and let
$\xi$ be a generator of $H_2(X)$. However, it is not true in general that this bound is the best possible for any
 rational homology $\cp^2$. As an example, consider Mumford surface (cf.~\cite{mumford:cp2}), which is an
algebraic surface of general type with the canonical class $3$. According to the generalized symplectic Thom
conjecture (see \cite{mst}) the minimal genus in the class of multiplicity $d>0$ equals
$$\frac{(d+1)(d+2)}{2}.$$
\end{remark}


\section{The main theorem}
In the previous section we saw how the moduli space of \sws\ on the complement of an embedded surface,
representing a given homology class in $\cp^2$, can be used to derive a lower bound for the genus of any smooth
representative of this class. Below we generalize this result to a bigger class of four-manifolds with $b_2^+=1$.
The bound is of course not explicit, but it can be effectively computed in many specific cases. The idea of the
proof is analogous to Kronheimer's proof of Donaldson's Theorem on definite intersection forms of closed
manifolds (see \cite{salamon}).
\begin{theorem}\label{main-theorem}
Let $X$ be a smooth closed oriented four-manifold with $b_1(X)\allowbreak =0$ and $b_2^+(X)=1$. Suppose $\Sigma \subset X$ is a
smooth surface representing homology class $d\xi$, where $d>1$ and $\xi \in H_2(X)$ is a primitive class of
positive self-intersection. Let $K$ be the set of all characteristic numbers $k \in \{0,\ldots,d\xi^2\}$
for $(X,\xi,d)$ which satisfy the following condition: there exists a \spinc structure on $X$ such that
\begin{equation}\label{X-condition}
c_1(Det)^2> \sigma(X) +4kd,
\end{equation}
and $\langle c_1(Det),\xi\rangle =k+d\xi^2$, where $Det$ is the determinant line of the \spinc structure. Suppose
that $K$ is not empty and let $k_0$ be the maximum of $K$. Then the genus $g$ of $\Sigma$ satisfies
\begin{equation}\label{genus-bound}
g>\frac{k_0d}{2}.
\end{equation}
\end{theorem}
\begin{proof}
We will use our standard notation $Z=X-\Sigma$ for the cylindrical end ma\-ni\-fold and $Y$ for the boundary of a
tubular neighborhood (oriented as the `boundary' of $Z$). Choose a regular perturbation $\omega \in \Omega$ (see
Proposition \ref{smoothness}); then the irreducible part $\MM^*$ of the perturbed moduli space $\MM$ is smooth.
Notice that (\ref{X-condition}) is equivalent to the dimension of the (perturbed) based moduli space $\widetilde\MM$
being positive (cf. Corollary \ref{formal-dim}); since this dimension is even, this also
implies that the dimension of $\MM$ is positive. With this remark, the statement of the theorem is equivalent to the
following claim, which we prove below: if the moduli space of \sws\ on $Z$ is positive dimensional for a given \spinc
structure, then it is not compact.

Suppose contrary to the statement of the theorem that for some $k \in K$, $kd/2 \ge g$. Fix a \spinc structure on
$X$ for which (\ref{X-condition}) holds for $k$. According to part (b) of Theorem \ref{3-dim-sol} there are no
irreducible solutions on $Y$ in the induced \spinc structure and therefore the moduli space $\MM$ is
compact (cf. Proposition \ref{compactness}).

Next we show that $\MM$ contains a unique reducible point $[A,0]$. Since $\omega\in\Omega$, $d^+$ is a
surjection from the space of extended one-forms to the space of self-dual two-forms. Hence the equation $F_A^+=\omega$
has a solution. Suppose now that $(A',0)$ is another solution and write $A'=A+\i\alpha$
for some one-form $\alpha$ on $Z$; clearly $d^+\alpha=0$ and therefore $d\alpha=0$. By choice of gauge the
asymptotic value $h$ of $\alpha$ is harmonic. Since the class
of $\alpha$ is trivial in $H^1(Z;\R)$, $h$ represents the trivial class in $H^1(Y;\R)$, hence $h=0$. Thus
there exists a function $f$ on $Z$, exponentially decaying to 0 at infinity, such that $\alpha=2df$, i.e.,
$(A,0)$ and $(A',0)$ are gauge equivalent.

To finish the proof of the claim, we need to understand the structure of $\MM$ at the reducible point $[A,0]$. Recall
that the index of the ASD part of the linearization at $(A,0)$ is zero. Since its kernel also vanishes (by an
argument similar to the one in the previous paragraph), the (Zariski) tangent space and the obstruction space at
$(A,0)$ correspond to the
kernel and cokernel of $D_A$ respectively. As we assumed that the index is positive, Proposition
\ref{smoothness} implies that the cokernel of the Dirac operator vanishes, so
the based moduli space $\widetilde\MM$ is smooth. The action of the group $S^1$ of constant gauge transformations
on the kernel of $D_A$ is by complex multiplication, hence a closed neighborhood $\VV$ of $[A,0]$ in $\MM$ is a cone
on some projective space $\cp^n$. Let $\NN$ be the smooth compact submanifold of $\MM$ with the boundary $\cp^n$,
obtained as the closure of the complement of $\VV$.  Denote by $c$ the Chern form of the $S^1$ bundle $\widetilde\NN
\to \NN$, where $\widetilde\NN \subset \widetilde\MM$ is the preimage of $\NN$. Note that the induced $S^1$ bundle
over the boundary $\cp^n$ is the tautological bundle.  So $\int_{\cp^n} c^n= \int_{\NN} d(c^n)=0$ is a contradiction.

We remark that for $k=0$, the condition $g>0$ that we obtain from the argument above is consistent with the
assumption that $\Sigma$ is  a sphere (needed for non-degeneracy of solutions on $Y$).

Finally, we check that it suffices to consider characteristic numbers $\le ad$, i.e., that for $ad<k<2ad$ we do not
get any new restrictions on the genus (here $a=\xi^2$). Suppose that for some $k$ in $(ad,2ad)$  condition
(\ref{X-condition}) holds
for some \spinc structure on $X$. Then we claim that there exists a \spinc structure on $X$ with the
characteristic number $k'=2ad-k$, for which (\ref{X-condition}) holds as well, hence $k' \in K$. To see this first
change the given \spinc structure by the line bundle $E$ on $X$ with $c_1(E)=-d\alpha$, where $\alpha$ is the
Poincar\'e dual of $\xi$. The characteristic number of the inverse of thus obtained \spinc structure is $k'$.  The
expression for the dimension of the moduli space is unaffected by these changes of the \spinc structure. For the first
change this follows from Corollary \ref{formal-dim}. For the second  note that $c_1(-Det)=-c_1(Det)$, hence the
class $c$, defined by  $\langle c_1(Det),\xi\rangle =\frac{p}{a}\alpha+c$, changes sign. The claim now follows from
the second form of the dimension formula (\ref{formal-dimension}). Since $k' \in K$, we have seen above that $g>
k'd/2=|k-2ad|d/2$; however, by part (b) of Theorem \ref{3-dim-sol}, this implies that there exist irreducible \sws\
on $Y$ in the \spinc structure induced by the one given on $X$ and hence the moduli space need not be compact.
\end{proof}
\begin{remark}
We note that the leading term in the genus bound for a divisible class $d\xi$, obtained from the above theorem,
equals $(d\xi)^2/2$, which is by a factor of 2 better than the bounds obtained via the
$g$-signature Theorem (cf. Rohlin  \cite{rohlin:genus}).
\end{remark}

A special case of interest occurs when the class $\xi \in H_2(X)$ is {\em characteristic}, that is, its
Poincar\'e dual is a characteristic class.
\begin{corollary}\label{characteristic-bound}
With notation as in Theorem \ref{main-theorem}, assume that $H_1(X)\allowbreak =0$, the signature of $X$ is negative,
and that $\xi$ is characteristic. Then the genus $g$ of any smooth surface $\Sigma$ representing $d\xi$ satisfies
$$g>\frac{[\Sigma]^2}{2}-\frac{d\xi^2}{2}=\binom{d}{2}\xi^2.$$
\end{corollary}
\begin{proof}
By Furuta's 10/8 Theorem \cite{furuta:11/8}, $X$ is odd.
Consider the \spinc structure on $X$ characterized by $c_1(Det)=(2d-1)\alpha$, where $\alpha$ is the Poincar\'e dual
of $\xi$. From $\langle c_1(Det),\xi\rangle = (2d-1)a=k+ad$ we get $k=a(d-1)$, where $a=\xi^2$.
Since the class $c$ is equal to 0 in
this case, (\ref{X-condition}) is equivalent to $(k-ad)^2>a\sigma(X)$, which is clearly true for
$k=ad-a$ and $\sigma(X) < 0$. This implies $g>a(d-1)d/2$.
\end{proof}


\section{Geometric constructions}\label{constr}
Let $X$ be a smooth four-manifold. For a class $\xi \in H_2(X)$ denote by $g_{\xi}(d)$ the minimal genus of a
smooth representative of $d\xi$; we write $g_{\xi}$ for $g_{\xi}(1)$. In this section we will show, using some simple
geometric constructions, that asymptotically $g_{\xi}(d)$ does not grow faster than $(d\xi)^2/2$. Combining this
with genus bounds from  Theorem \ref{main-theorem} we conclude that $(d\xi)^2/2$ describes the dominant
term in $g_{\xi}(d)$ (as $d \to \infty$) in a manifold $X$ with $b_2^+(X)=1$.
\begin{proposition}\label{construction}
Let $X$ be a smooth four-manifold and $\xi \in H_2(X)$ a class of positive self-intersection. Then
\begin{equation}\label{upper-bound}
g_{\xi}(d) \le \frac{(d\xi)^2}{2}-\bigg(\frac{\xi^2}{2}+1-g_{\xi}\bigg)d+1
\end{equation}
for any $d>1$. Moreover, there exists a smooth representative of $d\xi$ with the genus given by the right-hand side in
the above inequality.
\end{proposition}
\begin{proof}
Let $\Sigma \subset X$ be a smooth embedded surface of genus $g_{\xi}$ representing $\xi$ and let $\Sigma'$ denote
$\Sigma$ with the interiors of $a:=\xi^2$ disjoint disks removed. Think of the normal bundle $\nu_\Sigma$ of $\Sigma$
in $X$
as being obtained from the product bundle over $\Sigma'$ by adding a degree 1 bundle over a 2-disk for each puncture.
To construct $d$ copies $\Sigma_i$ of $\Sigma$ in general position, we choose $d$ distinct parallel copies $\Sigma_i'$
of $\Sigma'$. Over each 2-disk we cap-off $\Sigma_i'$ by adding a disk in such a way that any two disks intersect
transversely in a single point and any intersection point is common to two disks only. It is clear from the
construction that the total number of intersection points between the surfaces $\Sigma_i$ thus obtained equals
$a\binom{d}{2}$. A neighborhood of each intersection point looks like a pair of transversely intersecting disks.
Removing these and replacing them by annuli gives $\Sigma(d)$. Since the surfaces, obtained from $\Sigma_i$,
$i=1,\ldots,d$ by removing small (disjoint) disks around the intersection points are disjoint, we need $d-1$ annuli to
make a connected surface; each of the remaining annuli increases the genus by 1.
\end{proof}

As an immediate consequence of the above inequality we obtain the following bound on the genus of a representative
of a primitive class.
\begin{corollary}\label{primitive-bound}
Let $X$ be a smooth four-manifold and let $\xi \in H_2(X)$ be a class of positive self-intersection. If
$g_\xi(d)>(d\xi^2-\Delta_d)d/2$ for some $d>1$, then $g_\xi>(\xi^2-\Delta_d)/2$.
\end{corollary}

Assume now that $X$ is a smooth closed oriented four-manifold with $b_1(X)\allowbreak
=0$ and $b_2^+(X)=1$. If $\xi \in H_2(X)$
is a primitive class of positive self-intersection, then for any integer $d>1$ we know from Theorem
\ref{main-theorem} that $g_\xi(d)>(ad-\Delta_d)d/2$ for some $\Delta_d$. For example, if the signature of $X$ is
negative and  $\xi$ is characteristic, we can take $\Delta_d=a$ (see Corollary \ref{characteristic-bound}).
The previous corollary then implies that a characteristic class in $X$ is not represented by an embedded sphere.

In general we obtain the following consequence of Theorem \ref{main-theorem}. This result is equivalent to
the {\em generalized adjunction inequality} of Theorem \ref{adjunction-inequality}.
\begin{corollary}\label{form-of-bound}
Let $X$ be a smooth closed oriented four-manifold with $b_1(X)=0$ and $b_2^+(X)=1$. For any (primitive) class
$\xi \in H_2(X)$  of positive self-intersection  there exists a $\Delta \ge 0$ so that
$$g_\xi(d) > \frac{(d\xi^2-\Delta)d}{2}$$
for all $d \ge 1$.
\end{corollary}
\begin{proof}
We first check that the set $K$ in the statement of Theorem \ref{main-theorem} is non-empty for $d$ large enough.
Choose
a \spinc structure on $X$ with $c_1(Det)=2d\alpha-\gamma$, where $\gamma$ is a characteristic vector satisfying
$\gamma^2>\sigma(X)$ and $\gamma \cup \alpha \ge 0$. Such a characteristic vector clearly exists -- starting with
any characteristic vector we can get one satisfying these conditions by adding to it a large enough multiple of
$\alpha$. Then $k=ad-\langle\gamma ,\xi\rangle$ belongs to $[0,ad]$ for large enough $d$. We need to check that
(\ref{X-condition}) also holds:
$$c_1(Det)^2-\sigma(X)-4kd=\gamma^2-\sigma(X)>0.$$

Fix some $d$ for which $K$ is non-empty and a characteristic number $k \in K$ for $(X,\xi,d)$; let $\Delta:=ad-k$.
Denote by
$c_1=c_1(Det)$ the Chern class of the \spinc structure that satisfies (\ref{X-condition}) with $k$ and $d$, and let
$c_1(Det')=c_1+2n\alpha$ for some integer $n$. Then $\Delta'=a(d+n)-k'=\Delta$ and it follows from the second form of
the dimension formula (\ref{formal-dimension}) that the \spinc structure with determinant $Det'$ satisfies
(\ref{X-condition}) with $k'$ and $d+n$. This implies that the genus bound in the statement of the corollary holds
for all multiplicities $d+n$ for which $k'\ge 0$; it clearly holds for the rest.
\end{proof}


\section{Manifolds with signature zero}\label{signature=0}
In this section $X$ denotes a smooth closed oriented four-manifold with $b_1(X)=0$, $b_2^+(X)=1$ and signature
$\sigma(X)=0$.  Up to isomorphism, there are only two possible intersection forms such a manifold can have,
distinguished  by the parity. The even intersection form is given by $H=\left[\begin{matrix}0 & 1 \\ 1 &
0\end{matrix}\right]$ and is realized for example by $S^2 \times S^2$; the odd intersection form is given by
$E=\left[\begin{matrix}1 & 0 \\ 0 & -1 \end{matrix} \right]$ and is realized for example by $\cp^2 \# \cpbar$.
Genus bounds for divisible classes in $X$ that follow from Theorem \ref{adjunction-inequality} depend only on the
intersection pairing of the manifold, so we only need to consider two cases.
\begin{proposition}\label{bounds-for-H}
Suppose the intersection pairing of $X$ is isomorphic to $H$; let $\{\xi_1,\xi_2\}$ be a basis of $H_2(X)$ modulo
the torsion subgroup with respect to which the intersection pairing is given by $H$. Then any class $\xi \in H_2(X)$,
whose image in $H_2(X;\R)$ is given by $p\xi_1+ q\xi_2$ with $pq\ne 0$, satisfies
$$g_\xi\ge (|p|-1)(|q|-1).$$
\end{proposition}
\begin{proof}
After possibly changing the orientation of $X$ we may assume that the self-intersection $2pq$ of $\xi$ is positive;
then we  can further assume $p,q>0$. Denote by $\xi_i^* \in H^2(X;\R)$ the hom-dual of $\xi_i$ and consider a
characteristic vector $c$ satisfying $c=2\xi_1^*+2\xi_2^*$ in $H^2(X;\R)$. Since $c^2>0$, the claimed bound
follows from adjunction inequality (\ref{adjunction}).
\end{proof}

The following result for $S^2 \times S^2$ and the corresponding result regarding classes in $\cp^2 \# \cpbar$
(see below) was proved independently by Ruberman \cite{ruberman:thomrat} and by Li-Li \cite{li:mingenus}.
\begin{corollary}
With notation as in the previous proposition, assume further that $\xi_1$ and $\xi_2$ are represented by spheres
which intersect transversely in one point. Then
$$g_\xi= \left\{\begin{array}{ll}
(|p|-1)(|q|-1) & \hbox{ if }\ \ pq \ne 0 \\
0 & \hbox{ if }\ \ pq=0
\end{array} \right. .$$
\end{corollary}
\begin{proof}
Let $\Sigma_i$ be a representative of $\xi_i$ as in the statement.
Since the self-intersection of $\xi_i$ is zero, any class with $pq=0$ is represented by a sphere. Suppose now that
$pq \ne 0$; we may assume $p,q>0$. To construct a representative of $\xi$ with genus $(p-1)(q-1)$, take $p$ disjoint
copies of $\Sigma_1$ and $q$ disjoint copies of $\Sigma_2$, so that any copy of $\Sigma_1$ intersects any copy of
$\Sigma_2$ in exactly one point. Resolving the intersection points gives the required representative.
\end{proof}

\begin{proposition}\label{bounds-for-E}
Suppose the intersection pairing of $X$ is isomorphic to $E$; let $\{\xi_1,\xi_2\}$ be a basis of $H_2(X)$ modulo
the torsion subgroup with respect to which the intersection pairing is given by $E$. Then any class $\xi \in H_2(X)$
of positive self-intersection satisfies
$$g_\xi > \frac{p^2-q^2-3|p|+|q|}{2},$$
where $p\xi_1+ q\xi_2$ is the image of $\xi$ in $H_2(X;\R)$.
\end{proposition}
\begin{proof}
We may assume (by possibly changing the sign of $\xi_i$) that $p > q \ge 0$. Let $c$ be a characteristic class whose
real image is $3\xi_1^*-\xi_2^*$, where $\xi_i^* \in H^2(X;\R)$ denotes the hom-dual of $\xi_i$. As $c^2>0$, the
adjunction inequality implies the claimed genus bound, except for $q=0$ and $p \le 2$. In the latter cases the
claimed bound states $g_\xi \ge 0$, which is the best possible bound since in $\cp^2 \# \cpbar$ these classes are
represented by spheres.
\end{proof}
\begin{corollary}
With notation as in the previous proposition, assume further that $\xi_i$ are represented by disjoint spheres. Then
$$g_\xi= \left\{\begin{array}{ll}
\displaystyle{\frac{p^2-q^2-3|p|+|q|}{2}+1} & \hbox{ if }\ \  |p|>|q|  \\
\displaystyle{\frac{q^2-p^2-3|q|+|p|}{2}+1} & \hbox{ if } \ \ |q|>|p|  \\
0 & \hbox{ if }\ \  |p|=|q|
\end{array} \right. .$$
\end{corollary}
\begin{proof}
Let $\Sigma_i$ be a representative of $\xi_i$ as in the statement.
Note that any class of the form $(\pm 1,\pm 1)$ is represented by a sphere of self-intersection 0. Hence any class
$(p,q)$ with $ |p|=|q|$ is represented by a sphere.

Suppose that $ |p|>|q|$ (the remaining case is analogous). We may assume $p>q>0$ for the purpose of construction;
if $q=0$, the situation is as in $\cp^2$. To construct a representative of $\xi$ with stated genus, decompose $\xi$ as
$(p,q)=q(1,1)+(p-q)(1,0)$. Represent $q(1,1)$ by $q$ disjoint spheres, and $(p-q)(1,0)$ by a surface $\Sigma$ of genus
$(p-q-1)(p-q-2)/2$ which intersects each of the spheres in $p-q$ points. Finally resolve the intersection points.
\end{proof}


\section{Manifolds with negative signature}
Let $X$ be a smooth closed oriented four-manifold with $b_1(X)=0$, $b_2^+(X)=1$ and signature $\sigma(X)=1-n$
with $n \ge 2$. We will assume that the intersection form of $X$ is odd. This is true for $n \le 8$, since any such
form is odd. Without the restriction on $n$, the assumption holds for manifolds without 2-torsion according to Furuta's
$10/8$ Theorem  \cite{furuta:11/8}.

Fix a primitive class $\xi \in H_2(X)$ of positive self-intersection and choose a basis $\{\xi_0,\ldots,\xi_n\}$ of
$H_2(X)$ (modulo the torsion) with respect to which the intersection form is given by $\langle 1 \rangle \oplus n
\langle -1 \rangle$, and $\xi=(p,q_1,\ldots,q_n)$ with $p>0$ and $q_i \ge q_{i+1}\ge 0$; then $\xi^2=p^2-\sum
q_i^2$. Denote by $m=m_\xi$ the number of non-zero $q_i$'s.

It turns out that the genus bounds obtained from adjunction inequality (\ref{adjunction}) with $c$ the `canonical
class' $(3,-1,\ldots,-1)$ are only optimal for reduced classes. The notion of a reduced class, used by Li-Li
\cite{li:mingenus} to study genus bounds in rational surfaces $\cp^2 \# n\cpbar$ for $n\le 9$, extends naturally to
manifolds we are considering.
\begin{defn}
A class $\xi \in H_2(X)$ as above is called {\em reduced} with respect to the basis  $\{\xi_0,\ldots,\xi_n\}$
provided $m_\xi \le 9$ and $p \ge q_1+q_2 +q_3$, where $q_3=0$ if $n=2$.
\end{defn}
It is proved in \cite[Lemma 4.1]{li:mingenus} that in rational surfaces with $n \le 9$ any class of positive
self-intersection can be mapped to a reduced class (with respect to the standard basis) by an orientation preserving
diffeomorphism. The argument there also proves that any class $\xi \in H_2(X)$ with $m_\xi \le 9$ is reduced with
respect to some basis as above.

\begin{proposition}\label{neg-signature}
With above notation, suppose that $2 \le m \le 9$ and $\xi$ is reduced. Then
$$g_{\xi}(d)>\frac{(d\xi^2-3p+\sum q_i)d}{2}.$$
Moreover, $g_{\xi}(d)>0$ unless $d=1$, $m=2$ and $\xi=(p,p-1,1)$ for some $p>1$. Excluding the latter classes, given
any $g>0$ there is only a finite number of reduced classes with minimal genus no greater than $g$.
\end{proposition}
\begin{remark}
The last statement gives an affirmative answer to a conjecture of Li and Li \cite{li:mingenus}. In fact, B.H.~Li
proved in \cite{li:mingenreps} that the above lower bound is sharp in rational surfaces $X=\cp^2 \# n\cpbar$ for
$n \le 9$.
\end{remark}
\begin{proof}
Let $c$ be a characteristic vector whose real image with respect to the hom-dual basis is $(3,-1,\ldots,-1)$; clearly
$c^2>\sigma(X)$. Let $\Delta=\langle c,\xi \rangle=3p-\sum q_i$, where the sum, as all other sums over $i$,
runs from 1 to $m$. Since $\xi$ is reduced, $\Delta\ge 0$.

Suppose first that $m \ge 3$. Since $\xi$ is reduced,
$$p-\frac{3}{2} \ge (q_1- \frac{1}{2}) + (q_2- \frac{1}{2}) + (q_3- \frac{1}{2}) .$$
Using this along with $q_i \ge q_{i+1}$, we obtain
$$\xi^2-\Delta \ge (9-m)(q_3^2-q_3) \ge 0,$$
thus $g_{\xi}(d)\ge 1$ for any $d \ge 1$. Note in general that to establish the last claim of this proposition, it
suffices to show that an upper bound on minimal genus implies an upper bound on $q_1$; this is enough since
increasing the value of $p$ (while keeping $q_i$'s fixed) increases $\xi^2-\Delta$. The argument is simple if $m \le 8$,
as then $\xi^2-\Delta >(q_1-1/2)(q_2-1/2)$. The last inequality implies that an upper bound on minimal genus yields an upper
bound on $q_1$. For $m=9$ the first inequality of this paragraph gives
$$\xi^2-\Delta \ge \sum_{i=4}^7 (q_1-q_i)(q_i-\frac{1}{2})+\sum_{i=8}^9 (q_2-q_i)(q_i-\frac{1}{2});$$
we need to consider several cases. If $q_1>q_7$, a bound on minimal genus implies a bound on $q_1$. Same holds if
$q_1=q_7$, but $q_1>q_9$. Finally, if $q_1=q_9=q$,
$$\xi^2-\Delta =(p-3q)(p+3q-3),$$
but positive square condition implies $p>3q$; again it follows that there is only a finite number of such vectors
whose minimal genus is at most $g$.

If $m=2$, then $p \ge q_1+q_2$ implies
$$\xi^2-\Delta +2 \ge  2(q_1-1)(q_2-1),$$
which is strictly positive unless $q_2=1$. Note also that for $q_2>1$ there are only finitely many classes $(p,q_1,q_2)$
with minimal genus at most $g$. For $q_2=1$, we get
$$\xi^2-\Delta+2 = (p-q_1-1)(p+q_1-2);$$
this equals 0 only if $p=q_1+1$, and for $p\ge q_1+2$ there are only finitely many classes $(p,q_1,1)$ with minimal
genus no greater than $g$.
\end{proof}

It is easy to prove that the bound in the above proposition is the best possible bound obtainable from Theorem
\ref{adjunction-inequality}. However, this also follows from work of Li-Li \cite{li:mingenus} and
\cite{li:symplgenus}, where they prove that the minimal genus bound for a reduced class $\xi=(p,q_1,\ldots,q_m)$ of
positive self-intersection in a rational surface $\cp^2 \# n\cpbar$ with $n\le 9$  is given by
$$g_{\xi}(d)=\frac{(d\xi^2-3p+\sum q_i)d}{2}+1.$$
In proposition below we give a construction of a minimal genus representative in a special case.
\begin{proposition}
Let $X$, $\xi_i$, $\xi=(p,q_0,\ldots,q_m)$ and $m \le 9$ be as above, and let $\bar m$ be the largest value of $i$ for
which $q_i>2$. Assume further that $\sum_{i=1}^{\bar m} q_i \le p$ and that $\xi_i$ for $i=0,\ldots,m$ are represented
by  disjoint spheres. Then $\xi$ has a representative of genus
$$g_{\xi}=(\xi^2-3p+\sum q_i)/2+1.$$
\end{proposition}
\begin{proof}
Assume first that $\bar m =m$ and let $q:=\sum_{i=1}^{\bar m} q_i$. Then $\xi$ can be decomposed as $(p-q)\xi_0+\sum
q_i(\xi_0+\xi_i)$. By assumption $\xi_0$ and $\xi_0+\xi_i$ (for $i=1,\ldots,m$) can be represented by spheres
$\Sigma_i$ any two of which intersect transversely in one point. Moreover, the spheres $\Sigma_i$ for $i \ge 1$ have
self-intersection zero. To construct a representative for $\xi$, take $p-q$ copies of $\Sigma_0$ and $q_i$ disjoint
copies of $\Sigma_i$, so that the whole collection of spheres is in general position and any two spheres that intersect
have exactly one point in common. Note that the total number of intersection points of these $p$ spheres is
$$\binom{p-q}{2}+q(p-q)+\sum_{i<j} q_iq_j,$$
so after resolving the intersection points we obtain a minimal genus representative.

If $m> \bar m$, for any $q_i=2$ take two spheres representing $\xi_i$ that intersect transversely in one point.
Then connect one of the two spheres representing $\xi_i$ to a sphere $\Sigma_j$ for some $j \le \bar m$, obtaining
a surface $\Sigma$. Now cancel the $-1$ intersection point with one of the $+1$ intersection points of $\Sigma_j$.
To this end choose a curve $\gamma \subset \Sigma$ connecting the two intersection points and replace the complements
of small disks around the intersection points (cut out from the other surfaces, not $\Sigma$) by a tube which is the
restriction to $\gamma$ of the normal circle bundle of $\Sigma$. Resolving the remaining intersection points again
gives a minimal genus representative. Finally, if $q_i$ is $1$, connect the corresponding
sphere to the surface constructed before.
\end{proof}


\section{Geometric intersections of surfaces}\label{configurations}
Let $X$ be a smooth closed connected four-manifold. We say that a collection
of classes $\xi_1,\ldots,\xi_n \in H_2(X)$ is {\em algebraically disjoint} if $\xi_i \cdot \xi_j=0$ for all pairs
$i \ne j$. Classes in an algebraically disjoint collection can clearly be represented by disjoint surfaces -- starting
with any choice of representatives in general position, we can eliminate a pair of $\pm 1$ intersection points between
two surfaces by adding a one-handle to one of them. An important question is whether the classes can be represented
by disjoint surfaces of low genus. It turns out that the minimal genus representatives of the classes intersect in
general and we will derive a lower bound for the number of pairs of $\pm 1$ intersection points.
\begin{proof}[Proof of Theorem  \ref{n-adj-ther}]
We assume for the purpose of the proof that $\Sigma_i$ is representing a divisible class $d_i\xi_i$ for some $d_i>1$.
The result for $d_i=1$ then follows as in Corollary \ref{primitive-bound}.
Consider $Z=X-\cup_{i=1}^n \Sigma_i$ as a cylindrical end manifold with cylindrical ends $[0,\infty) \times Y_i$,
where $Y_i$ is the boundary of a tubular neighborhood of $\Sigma_i$ (with the opposite orientation). Note that $Z$ is
negative semi-definite, so we can adapt the argument that we used to derive genus bounds to this context.

We work with exponentially decaying \sw\ configuration space on $Z$ with reducible asymptotic values; the
background connection on $Z$ agrees on the end $[0,\infty) \times Y_i$ with the pull-back of the product connection
on $Y_i$ (cf.~the proof of Proposition \ref{fiber-index}). By gauge fixing
at infinity we can assume that the asymptotic values differ by imaginary-valued harmonic one-forms on each end.
We use gauge group based at infinity on the first end $[0,\infty) \times Y_1$. To compute the index of the
deformation complex, we follow the proof of Proposition \ref{fiber-index}. The quotient of the deformation complex
by the fiber complex is
$$ 0 \to \oplus_{i=2}^n \R \to \oplus_{i=1}^n \HH(Y_i) \to 0 \to 0,$$
where $\HH(Y_i)$ denotes the space of imaginary-valued one-forms on $Y_i$. Hence the index of the deformation complex
is the index of the fiber complex plus $2\sum g_i -(n-1)$, where $g_i$ is the genus of $\Sigma_i$. The computation
of the fiber index splits in the anti self-dual part and the Dirac part. Note for the first that each end contributes
1 to the spectral flow; it follows that the fiber index of the ASD part equals $-2\sum g_i +n-1$. The index of the
deformation complex is thus given by
$$\frac{1}{4}(c_1(Det)-\sigma(X))-\sum_{i=1}^n k_id_i$$
(see also Lemma \ref{c1^2-flat}),
where $Det$ is the determinant line of a \spinc structure on $X$ satisfying
$$k_i=\langle c_1(Det),\xi_i \rangle - d_i\xi_i^2 \in [0,d_i\xi_i^2]$$
for all $i$.

Next we argue as in the proof of Theorem \ref{main-theorem}. The based moduli space of \sws\ on $Z$ is smooth (after
choosing an appropriate perturbation) and contains a unique reducible point. From the structure of the moduli
space close to the reducible point we conclude that if the moduli space is positive dimensional, then it is not
compact. Therefore $g_i>k_id_i/2$ for at least one $i$. Finally, setting $c_1(Det)=\sum 2d_i\alpha_i -c$, where
$\alpha_i$ is the Poincar\'e dual of $\xi_i$, we see that the formal dimension of the (based) moduli space is positive
if and only if $c^2>\sigma(X)$ (see the proof of Corollary \ref{form-of-bound}). The condition at the end of the
previous paragraph then becomes $\langle c,[\Sigma_i] \rangle \in [0,[\Sigma_i]^2]$, while the genus bound can be
expressed as $\chi(\Sigma_i)+[\Sigma_i]^2 \le \langle c,[\Sigma_i]\rangle$.
\end{proof}

In what follows we will restrict our attention to manifolds with $b_2^+(X)=2$ in order to keep the discussion simple.
\begin{theorem}\label{intersections}
Let $X$ be a smooth closed connected four-manifold with $b_1(X)=0$ and $b_2^+(X)=2$, and let  $\Sigma_1,
\Sigma_2$ be embedded surfaces in general position, representing algebraically disjoint classes
of positive self-intersection. Suppose that a characteristic vector $c \in H^2(X)$ satisfies
$$c^2>\sigma(X),\ \ \ \ \langle c,[\Sigma_i] \rangle \ge 0 \ \ \hbox{for $i=1,2$}, $$
$$\hbox{and} \ \ \ \chi(\Sigma_i)+[\Sigma_i]^2 > \langle c,[\Sigma_i]\rangle  \ \ \hbox{for $i=1,2$}.$$
Then
$$g(\Sigma_1)+g(\Sigma_2)+N \ge \frac{[\Sigma_1]^2+[\Sigma_2]^2- \langle c, [\Sigma_1]+[\Sigma_2] \rangle}{2}+1,$$
where $N$ denotes the number of pairs of $\pm 1$ intersection points between $\Sigma_1$ and $\Sigma_2$.
\end{theorem}
\begin{proof}
If condition (\ref{n-adjunction}) fails for both surfaces, then by Theorem \ref{n-adj-ther} the surfaces are not
disjoint.
To  construct disjoint representatives, we trade pairs of $\pm 1$ intersection points for one-handles --  this way
ellimination of a pair of intersection points increases the genus of one of the surfaces by 1. We add the maximal
possible number of handles to $\Sigma_1$, so that the resulting surface still does not satisfy (\ref{n-adjunction}),
and add the rest of the handles to $\Sigma_2$. Since the sum of the genera of thus constructed disjoint surfaces
equals $g(\Sigma_1)+g(\Sigma_2)+N$, the claimed inequality follows from Theorem \ref{n-adj-ther}.
\end{proof}

We compare this to bounds obtained using g-signature Theorem. As is the case for genus bounds, the result we obtained
is roughly by a factor of 2 better. Specifically, we state the following consequence of a Theorem of Gilmer
\cite{gilmer:thesis}.
\begin{proposition} {\bf (Gilmer \cite{gilmer:thesis})}
Let $X$ be a smooth closed connected four-manifold with $b_1(X)=0$ and $b_2^+(X)=2$, and let  $\Sigma_1,\ \Sigma_2$ be
embedded surfaces in general position, representing algebraically disjoint classes of positive
self-intersection. If $\Sigma_1$ and $\Sigma_2$ are not disjoint and $[\Sigma_1]+[\Sigma_2]$ is divisible by $2$, then
$$g(\Sigma_1)+g(\Sigma_2)+N \ge \frac{[\Sigma_1]^2+[\Sigma_2]^2}{4}-1,$$
where $N$ denotes the number of pairs of $\pm 1$ intersection points between $\Sigma_1$ and $\Sigma_2$.
\end{proposition}

\subsection{Examples}
Let $X=S^2 \times S^2 \# S^2 \times S^2$ and let $\xi_1=(p,q,0,0)$ and $\xi_2=(0,0,r,s)$ be classes of
positive self-intersection, expressed with respect to the standard basis for $H_2(X)$. We may assume that
 $p,q,r,s> 0$. If $\xi_i$ is primitive, it is represented by an embedded sphere in $X$,
according to a Theorem of Wall \cite{wall:diffeos}; however, for $p,q,r,s\ge 2$ it is not represented by a sphere in its
summand.
Let $\Sigma_1$ and $\Sigma_2$ be smooth representatives of $\xi_1$ and $\xi_2$ in general position. Denote by $g_i$
the genus of $\Sigma_i$ and  by $N$ the number of pairs of $\pm 1$ intersection points between $\Sigma_1$ and
$\Sigma_2$. Using characteristic vectors $c=(2,2,0,0)$ and $c=(0,0,2,2)$ in Theorem \ref{intersections} gives the
following lower bounds for $g_1+g_2+N$:
$$\begin{array}{lll}
(p-1)(q-1)+rs & \hbox{\ \ if\ \ } & p,q\ge 2,\  g_1<(p-1)(q-1) \hbox{\ \ and\ \ } g_2\le rs, \cr
pq+(r-1)(s-1) & \hbox{\ \ if\ \ } & r,s\ge 2,\  g_1\le pq \hbox{\ \ and\ \ } g_2<(r-1)(s-1).
\end{array}$$
In particular, if $\xi_1$ and $\xi_2$ with $p,q,r,s\ge 2$ are represented by spheres $\Sigma_1$ and $\Sigma_2$ in $X$,
then
$$N \ge \max \{(p-1)(q-1)+rs,pq+(r-1)(s-1)\}.$$

Consider now $X=\cp^2 \# \cp^2$ and let $\xi_1=(p,q)$ and $\xi_2=(q,-p)$ for some $p,q>0$, expressed with respect to
the standard basis for $H_2(X)$. Note that $\xi_i$ has a smooth representative $\Sigma_i$ of genus
$$g=\frac{p^2+q^2-3(p+q)}{2}+2,$$
obtained as the connected sum of minimal genus representatives for classes of divisibility $p$ and $q$ in $\cp^2$.
Let $\Sigma_1$ and $\Sigma_2$ be genus $g$ representatives of $\xi_1$ and $\xi_2$ in general position, and let $N$ be
the number of pairs of $\pm 1$ intersection points. Using characteristic vectors $c=(3,-1)$, $c=(3,1)$ and $c=(1,-3)$
in Theorem \ref{intersections}, we obtain
$$N \ge \left\{\begin{array}{lll}
p+2q-3 & \hbox{\ \ if\ \ } & p\ge 2 \cr
2p+q-3 & \hbox{\ \ if\ \ } & q\ge 2 \hbox{\ \ and\ \ } p \le 3q \cr
p+4q-3 & \hbox{\ \ if\ \ } & q\ge 2 \hbox{\ \ and\ \ } p \ge 3q
\end{array}\right. .$$
In particular, for $q=1$ we have $N \ge p-1$ and this bound is sharp, which can be seen as follows. Decompose
$(p,1)=(p-1)(1,0)+(1,1)$ and $(1,-p)=(1,-1)-(p-1)(0,1)$. Since the classes $(1,1)$ and $(1,-1)$ can be represented
by disjoint embedded spheres, $N=p-1$.

As the last example consider $X=2\cp^2 \# 2\cpbar$, and let $\xi_1=(p,0,q,0)$ for some $p>q>0$ and $\xi_2=(0,r,0,s)$
for some $r>s>0$. These classes are represented by spheres in $X$ (according to the Theorem of Wall mentioned above),
but not in their copy of $\cp^2 \# \cpbar$, unless $p=q+1$ or $r=s+1$. By choosing representatives $\Sigma_i$ for $\xi_i$
of small genus and in general position, we get the following bound based on Theorem \ref{intersections}
(using $c=(3,1,-1,-1)$ and $c=(1,3,-1,-1)$):
$$g_1+g_2+N \ge \left\{\begin{array}{lll}
\dsp{\frac{p^2-q^2-3p+q}{2}}+\dsp{\frac{r^2-s^2-r+s}{2}}+1 & \hbox{\ \ if\ \ } & p\ge q+2 \cr
\dsp{\frac{p^2-q^2-p+q}{2}}+\dsp{\frac{r^2-s^2-3r+s}{2}}+1 & \hbox{\ \ if\ \ } & r\ge s+2 \cr
\end{array}\right. .$$
Small genus here means that for a formula to hold, $g_1$ has to be no greater than the first summand and $g_2$ has
to be no greater than the second summand.

\begin{remark}
It is an interesting question whether the bounds obtained in the above examples are optimal.
\end{remark}



\providecommand{\bysame}{\leavevmode\hbox to3em{\hrulefill}\thinspace}

\end{document}